\documentclass[amssymb,11pt]{amsart}
\usepackage{amscd,amssymb,euscript,amsthm}
\theoremstyle{plain}
\newtheorem{thm}{Theorem}[section] 
\newtheorem{cor}[thm]{Corollary}
\newtheorem{prop}[thm]{Proposition}

\theoremstyle{definition}
\newtheorem{defn}[thm]{Definition}
\theoremstyle{remark}
\newtheorem{rem}[thm]{Remark}

\numberwithin{equation}{section}%assigns eqns section numbers
\newcommand{\ran}{\operatorname{ran}}

\newcommand{\tr}{\operatorname{trace}}

\def\<{\left<}
\def\>{\right>}
\def\cstar{$C^*$-algebra}
\begin{document}
\title[Standard Hilbert Modules]{Quotients of Standard Hilbert Modules}
\author{William Arveson}
\thanks{supported by 
NSF grant DMS-0100487} 
\address{Department of Mathematics,
University of California, Berkeley, CA 94720}
\email{arveson@math.berkeley.edu}
\begin{abstract}
We initiate a study of Hilbert modules 
over the polynomial algebra $\mathcal A=\mathbb C[z_1,\dots,z_d]$ 
that are obtained by completing $\mathcal A$ with respect 
to an inner product having certain 
natural properties.  
A {\em standard Hilbert module} is 
a finite multiplicity version 
of one of these.    Standard Hilbert 
modules occupy a position analogous to that of 
free 
modules of finite rank in commutative algebra, and  
their quotients by submodules give rise to universal solutions of 
nonlinear relations. 
Essentially all of the basic Hilbert modules 
that have received attention over the years 
are standard - 
including the Hilbert module of the $d$-shift, the Hardy and Bergman 
modules of the unit ball, modules associated with 
more general domains 
in $\mathbb C^d$, and those associated with 
projective algebraic varieties.

We address the general problem of determining when 
a quotient $H/M$ of an 
essentially normal standard 
Hilbert module $H$ is essentially normal.  This 
problem has been resistant.  Our main result 
is that it can be ``linearized" in 
that the nonlinear relations defining 
the submodule $M$ can be reduced, appropriately, to linear relations 
through an iteration procedure, 
and we give a concrete description 
of linearized quotients.   
\end{abstract} 
%
%\date{13 March, 2004}
%
\maketitle

\section{Introduction}\label{S:in}

Let $T_1,\dots,T_d$ be a commuting $d$-tuple of 
operators on a Hilbert space $H$.  One can 
view $H$ as a module over the polynomial algebra 
$\mathcal A=\mathbb C[z_1,\dots,z_d]$ in the 
natural way
$$
f\cdot\xi=f(T_1,\dots,T_d)\xi,\qquad f\in\mathcal A,\quad \xi\in H,   
$$
and such an $H$ is called a Hilbert module of 
dimension $d$, or simply a {\em Hilbert module} when the
dimension is understood.  A Hilbert module 
is said to be {\em essentially normal} if the self-commutators 
$T_k^*T_j-T_jT_k^*$ of its ambient operators are all compact, and 
more specifically,
{\em $p$-essentially normal} if the self-commutators belong 
to the Schatten class $\mathcal L^p$ - $p$ being a number in 
the interval $[1,\infty]$, where $\mathcal L^\infty$ is interpreted 
as the \cstar\ $\mathcal K$ of compact operators on $H$.  

As in algebra, examples in multivariable operator 
theory are constructed most naturally through the formation of 
quotients -- by passing from the simplest ``free" Hilbert modules 
$H$ to their quotient Hilbert modules $H/M$, in which $M$ is the  
closed submodule of $H$ generated by the relations one 
seeks to satisfy.  However, in the operator-theoretic 
context, properties like essential normality, Fredholmness, 
and finiteness of the cohomology of the associated Koszul complex,  
do not propagate naturally from $H$ to its submodules 
or quotients, with the result 
that it is unclear whether Hilbert modules with the sought-after 
properties exist.

For example, consider the problem of 
constructing operator solutions $X_1,\dots,X_d$ to systems of algebraic 
equations of the form 
\begin{equation}\label{inEq0}
f_k(X_1,\dots,X_d)=0,\qquad k=1,\dots,r,
\end{equation}
where the $f_k$ are homogeneous 
polynomials in $d$ variables and 
$X_1,\dots,X_d$ are mutually commuting bounded operators on a Hilbert space
$H$.  We can make the point well enough with the following 
special case, in which one fixes a positive integer $n\geq 2$ 
and seeks a commuting triple 
$X,Y,Z\in\mathcal B(H)$ with the property 
\begin{equation}\label{inEq1}
X^n + Y^n = Z^n.  
\end{equation}
Such a triple can be viewed as a ``noncommutative curve".  
We say {\em noncommutative} in order to rule 
out variants of classical solutions such as those 
obtained by assembling 
a sequence of numerical 
solutions $(x_k,y_k,z_k)\in\mathbb C^3$, $k\geq 1$, of 
$x^n+y^n=z^n$   
into diagonal operators acting on $\ell^2(\mathbb N)$ such as
$$
X=
\begin{pmatrix}
x_1&0&0&\cdots\\
0&x_2&0&\cdots\\
0&0&x_3&\cdots\\
\vdots&\vdots&\vdots
\end{pmatrix}
Y=
\begin{pmatrix}
y_1&0&0&\cdots\\
0&y_2&0&\cdots\\
0&0&y_3&\cdots\\
\vdots&\vdots&\vdots
\end{pmatrix} 
Z=
\begin{pmatrix}
z_1&0&0&\cdots\\
0&z_2&0&\cdots\\
0&0&z_3&\cdots\\
\vdots&\vdots&\vdots
\end{pmatrix}.  
$$

While in general, the operators 
$X,Y,Z$ of (\ref{inEq1}) must commute, 
the unital \cstar\ $C^*(X,Y,Z)$ that 
they generate is typically noncommutative.  
We 
seek solution triples $X,Y,Z$ that are ``universal" in 
an appropriate sense, 
that generate 
an irreducible \cstar\ containing the \cstar\ $\mathcal K$ of 
compact operators, and which is commutative modulo $\mathcal K$.  
The latter properties are summarized in an exact sequence of \cstar s 
\begin{equation}\label{inEq2}
0\longrightarrow \mathcal K\longrightarrow C^*(X,Y,Z)\longrightarrow 
C(V)\longrightarrow 0,   
\end{equation}
in which $V$ is a compact subset of  
$\{(x,y,z)\in\mathbb C^3: x^n+y^n=z^n\}$.  
The sequence (\ref{inEq2}) defines an extension of 
$\mathcal K$ by $C(V)$ in the sense of Brown-Douglas-Fillmore 
and, as shown in \cite{bdf}, it gives rise to an element 
of the $K$-homology of the space $V$.  Of course, 
similar comments 
apply to the more general 
``noncommutative varieties" defined by operator 
solutions $X_1,\dots,X_d$ of systems 
of equations of the form (\ref{inEq0}).

Following basic principles, 
one constructs solutions of algebraic equations such as 
(\ref{inEq1}) by passing from the polynomial algebra 
$\mathbb C[x,y,z]$ to its quotient by the ideal generated 
by $x^n+y^n-z^n$.  The elements $\dot x, \dot y, \dot z$ 
obtained by projecting $x,y,z$ to the quotient 
are forced to satisfy $\dot x^n+\dot y^n=\dot z^n$, 
and there is an obvious 
sense in which this solution $\dot x,\dot y,\dot z$ is the universal 
one.  

When one attempts to carry out this construction of solutions 
in the context of Hilbert space operators, one encounters 
an exasperating difficulty.  To illustrate the point, 
let us complete the polynomial 
algebra $\mathbb C[x,y,z]$ in an appropriate inner product --  
for definiteness we choose the inner product associated 
with the $3$-shift \cite{arvSubalgIII} - giving rise to a Hilbert 
space $H^2(x,y,z)$.  The operators $X,Y,Z$ on $H^2(x,y,z)$ defined 
by multiplication by the basic variables $x,y,z$ are bounded, 
we may form the closed invariant subspace $M\subseteq H^2(x,y,z)$ 
generated by $x^n+y^n-z^n$, and its Hilbert space quotient 
$H=H^2(x,y,z)/M$.     The operators 
$X,Y,Z$ promote naturally to operators 
$\dot X,\dot Y,\dot Z$ on the quotient $H$, the promoted
operators satisfy (\ref{inEq1}), 
and straightforward computation shows that the 
\cstar\ 
$C^*(\dot X,\dot Y,\dot Z)$ is irreducible and contains 
$\mathcal K$.  Thus, if one knew that all 
self-commutators $\dot X^*\dot X-\dot X\dot X^*$, 
$\dot X^*\dot Y-\dot Y\dot X^*$, $\dots$, 
were compact, then $C^*(\dot X,\dot Y,\dot Z)$ would 
be commutative modulo $\mathcal K$ and we would have 
an extension  of the desired form 
\begin{equation}\label{inEq3}
0\longrightarrow \mathcal K\longrightarrow C^*(\dot X,\dot Y,\dot Z)\longrightarrow 
C(V)\longrightarrow 0.  
\end{equation}
The  
difficulty is that {\em it is unknown if $C^*(\dot X,\dot Y,\dot Z)$ is 
commutative modulo $\mathcal K$ for any $n\geq 3$.}   

The difficulty only grows in more general settings.  
For example, suppose that $M$ is a graded submodule of the 
finite multiplicity Hilbert module 
$H^2(x,y,z)\otimes \mathbb C^r$, and let $H$ be the quotient 
Hilbert module 
$$
H=\frac{H^2(x,y,z)\otimes\mathbb C^r}{M}.
$$ 
This quotient $H$ no longer corresponds 
so simply to solutions of 
equations like (\ref{inEq1}), but could represent 
a Hilbert space 
of sections of a vector bundle or sheaf over an appropriate algebraic set.  
Interpretations aside, one still has the basic 
operator-theoretic question as 
to whether the \cstar\ $C^*(\dot X,\dot Y,\dot Z)$ 
generated by the natural operators of $H$ is commutative 
modulo $\mathcal K$.  Again, the answer is unknown 
in most cases of interest.  

The purpose of this paper is to initiate the study of a broad 
context in which one can confront this issue, 
and which contains all of the important examples.    
We show that the general  
problem of proving essential normality of quotient modules 
defined by nonlinear relations can be reduced to the case 
in which the relations that define the quotient module 
are {\em linear}.  Despite the explicit formulation 
of the linearized problem,  
it remains unsolved in general.  
Problems and conjectures are discussed in Section \ref{S:lq}.

Initially, we took up this program out of a desire to 
give a natural proof that 
the curvature invariant of \cite{arvCurv}, \cite{arvDirac} is stable under compact 
perturbations and homotopy, by showing that pure finite-rank 
graded $d$-contractions satisfy the Fredholm property.  In turn, 
that question led us to attempt to establish essential normality 
for certain quotient Hilbert modules (see Proposition \ref{pmProp1} below).  
The preceding discussion shows that other fundamental issues of operator 
theory lead naturally to the same problem.

\section{Graded Completions of $\mathbb C[z_1,\dots,z_d]$}\label{S:gc}

We consider Hilbert modules  
$G$ obtained by completing the algebra 
of complex polynomials $\mathcal A=\mathbb C[z_1,\dots,z_d]$ 
in an inner product with the property that the 
natural multiplication operators $Z_1,\dots,Z_d$ associated 
with the generators $z_1,\dots,z_d$ are bounded.  
One can increase multiplicity by 
forming the direct sum of $r<\infty$ copies 
$G\oplus \cdots\oplus G=G\otimes\mathbb C^r$ of 
such a module $G$, which is a Hilbert 
module whose natural operators are 
multiplicity $r$ versions of the original $Z_k$.  
It is convenient to abuse notation 
by also writing these finite multiplicity 
multiplication operators as $Z_1,\cdots,Z_d$.  

In this section we 
single out a class of inner products on the 
algebra $\mathcal A$ whose 
completions are Hilbert modules, the finite 
multiplicity versions of which 
form effective building blocks 
for multivariable operator theory in Hilbert 
spaces.  We call them {\em standard} Hilbert modules.  
This is a very broad class of Hilbert modules that  
differs in several ways from  
classes that have been previously studied 
\cite{mulVas}, \cite{dougMV1},  
\cite{dougMisEQ}.  
For example, 
standard Hilbert modules are not necessarily subnormal, 
nor are they 
necessarily associated with a reproducing kernel.  

Indeed, standard Hilbert modules occupy a position 
analogous to that of free modules in the 
algebraic theory 
of finitely generated modules over $\mathcal A$; 
they are basically the 
Hilbert modules that have 
the same {\em cohomology} as free modules in the 
algebraic theory (see Remark \ref{fpRem0.0}).  
However, while there is only 
one algebraic free module of rank one - namely $\mathcal A$ 
itself, there 
are many inequivalent standard Hilbert modules of rank one.  
That class of 
Hilbert modules includes all of the basic examples 
that have been studied in recent years, 
including the space $H^2$ of 
the $d$-shift, the Hardy space of the unit sphere 
in $\mathbb C^d$, the Bergman space of the unit ball, 
as well as Hilbert modules associated with 
other domains in $\mathbb C^d$ and 
projective algebraic varieties.

\begin{rem}[Graded inner products on $\mathcal A$]
There is a natural action $\Gamma$ of the circle group 
on $\mathcal A$ defined by 
$$
\Gamma(\lambda)f(z_1,\dots,z_d)=f(\lambda z_1,\dots,\lambda z_d), 
\qquad f\in\mathcal A,\quad \lambda\in\mathbb T.  
$$
We write $\mathcal A_n$ for the linear space of homogeneous 
polynomials of degree $n$
$$
\mathcal A_n=\{f\in\mathcal A: \Gamma(\lambda)f=\lambda^n f, 
\quad \lambda\in\mathbb T\},
\qquad n=0,1,2,\dots, 
$$
and the polynomial algebra decomposes into an algebraic 
direct sum of homogeneous subspaces 
$$
\mathcal A=\mathbb C\dotplus \mathcal A_1\dotplus \mathcal A_2\dotplus\cdots.  
$$

An inner product $\langle\cdot,\cdot\rangle$ on 
$\mathcal A$ is invariant under the action of $\Gamma$ 
$$
\langle \Gamma(\lambda)f,\Gamma(\lambda)g\rangle=\langle f,g\rangle,
\qquad f,g\in \mathcal A, \quad \lambda\in \mathbb T, 
$$
iff the  homogeneous spaces 
are mutually orthogonal: $\langle \mathcal A_m,\mathcal A_n\rangle=\{0\}$ for $m\neq n$.  
Such inner products are called {\em graded}.  
If, for 
each $k=1,\dots,d$, the 
multiplication operator $Z_k:f\mapsto z_kf$, $1\leq k\leq d$ 
is bounded relative to the norm $\|f\|=\langle f,f\rangle^{1/2}$, then the 
completion of $\mathcal A$ in this inner product is a 
Hilbert module over $\mathcal A$.  
\end{rem}

\begin{defn}\label{}
A {\em graded completion} of $\mathcal A$ is a Hilbert module 
$G$ obtained by completing $\mathcal A$ in a graded inner product 
with the additional property that the linear space 
$Z_1G+\cdots+Z_dG$ is closed.  
\end{defn}

Examples of graded completions include the Bergman and Hardy 
modules of the ball $\{z\in\mathbb C^d: \|z\|<1\}$, the space 
$H^2$ of the $d$-shift, the Bergman and Hardy modules of 
$d$-dimensional polydisks, and in fact most Bergman modules of 
more general domains in $\mathbb C^d$ that admit circular symmetry.  

\begin{rem}[Normalization of Coordinates]
In every graded completion $G$ of $\mathbb C[z_1,\dots,z_d]$, we 
can replace the indicated basis $z_1,\dots,z_d$ for $\mathcal A_1$ 
with another basis $\tilde z_1,\dots,\tilde z_d$, 
if necessary, to achieve the normalization 
$\langle \tilde z_i,\tilde z_j\rangle=\delta_{ij}$.  While this has the 
effect of changing the 
original set of multiplication 
operators $Z_1,\dots,Z_d$ into another linear basis 
$\tilde Z_1,\dots,\tilde Z_d$ for 
the operator space they span, it does 
not significantly affect properties of 
the Hilbert module $G$.  For 
example, one shows easily that if the self-commutators 
of the original operators $Z_j^*Z_k-Z_kZ_j^*$ belong to $\mathcal L^p$, 
then so do the self-commutators 
$\tilde Z_j^{*}\tilde Z_k-\tilde Z_k\tilde Z_j^*$.  Thus, we 
can assume throughout that 
{\em for any graded completion $G$, the coordinates $z_1,\dots,z_d$ are 
an orthonormal subset of $G$}.  
\end{rem}

For every graded completion $G$, the representation  
$\Gamma$ extends naturally to a 
strongly continuous unitary representation 
of $\mathbb T$ on $G$ (also written $\Gamma$) whose spectral subspaces 
$$
G_n=\{\xi\in G: \Gamma(\lambda)\xi=\lambda^n\xi\},\qquad n\in\mathbb Z
$$ 
vanish 
for negative $n$, satisfy $G_n=\mathcal A_n$ for $n\geq 0$, and one has 
$$
G=G_0\oplus G_1\oplus G_2\oplus\cdots.   
$$
The space $Z_1G+\cdots+Z_dG$ is the orthocomplement of the 
one-dimensional space of constants $\mathbb C\cdot1$.   
$\Gamma$ is called 
the {\em gauge group} of $G$; it relates to the ambient operators by way of 
$\Gamma(\lambda)Z_k\Gamma(\lambda)^*=\lambda Z_k$,  
$k=1,\dots,d$, $\lambda\in \mathbb T$.

\begin{rem}[Irreducibility]
A submodule $M\subseteq H$ of a Hilbert 
module $H$ is said to be {\em reducing} 
if its orthocomplement $M^\perp$ is 
also a submodule.  In this case we also 
refer to $M$ as a {\em summand} since 
it gives rise to a decomposition 
$H=M\oplus N$ of $H$ into a direct 
sum of Hilbert modules.  
$H$ is said to be 
{\em irreducible} if it 
has no nontrivial summands.  Equivalently, 
$H$ is irreducible iff the $*$-algebra 
generated by the ambient operators 
$Z_1,\dots,Z_d$ of $H$ has  
commutant $\mathbb C\cdot\mathbf 1_H$.  

\begin{prop}\label{gcProp1}
Every graded completion of $\mathcal A$  
is irreducible.  
\end{prop}

\begin{proof}
Let $G$ be a graded completion 
and let $P\in\mathcal B(G)$ be a projection that commutes with 
the coordinate operators $Z_1,\dots,Z_d$ 
of $G$.  Choose a unit vector $v\in G_0$, so that 
$G=[\mathcal Av]$.  It follows that 
$G_0=[v]$ is the orthocomplement of 
the closed subspace 
$$
Z_1G+\cdots+Z_dG.
$$  
Since $P$ commutes with each $Z_k$ it commutes with 
$Z_1Z_1^*+\cdots+Z_dZ_d^*$, and therefore with the projection 
onto the range of the latter operator, namely the subspace 
$Z_1G+\cdots+Z_dG$.  Hence $P$ commutes with 
the rank-one projection 
$v\otimes \bar v$, 
and we conclude that either $Pv=v$ or $Pv=0$. 
 
If $Pv=v$,  
then $P$ must restrict to the 
identity operator on 
the closed submodule $[\mathcal Av]=G$ 
generated by $v$, hence $P=\mathbf 1$.  
Similarly, $Pv=0$ implies $P=0$.  
\end{proof}
\end{rem}

\begin{rem}[Number operator]
The {\em number operator} of a graded completion $G$ is defined 
as the self-adjoint generator $N$ of the gauge group
$$
\Gamma(e^{it})=e^{itN},\qquad t\in\mathbb R.  
$$  
The number operator is self-adjoint, 
has integer spectrum 
$\{0,1,2,\dots\}$ and its minimal spectral projections are the 
projections onto the homogeneous spaces $G_n=\mathcal A_n$, 
$n=0,1,2,\dots$.  Since 
the dimensions of the spaces $\mathcal A_n$ do not 
depend on the inner product chosen, {\em any 
two graded completions of $\mathcal A$ have unitarily 
equivalent number operators}.  It was shown in the 
appendix of \cite{arvSubalgIII} that the number operator 
satisfies 
$$
(N+\mathbf 1)^{-1}\in\mathcal L^p \iff p>d.  
$$
\end{rem}

\begin{rem}[Graded Hilbert modules, Gauge groups]
More general graded Hilbert modules $H$ 
can be defined in two equivalent ways.  
One specifies 
either a $\mathbb Z$-grading 
for $H$ 
$$
H=\cdots\oplus H_{-1}\oplus H_0\oplus H_1\oplus\cdots
$$
in which $Z_kH_n\subseteq H_{n+1}$, for all $1\leq k\leq d$, 
$n\in\mathbb Z$, or one specifies a gauge group - 
a strongly continuous 
unitary representation $\Gamma$ of the circle 
group $\mathbb T$ on $H$ such that 
$$
\Gamma(\lambda)Z_k\Gamma(\lambda)^*=\lambda Z_k
\qquad 1\leq k\leq d,\quad \lambda\in\mathbb T.  
$$  
One passes back and forth via the associations 
$H_n=\{\xi\in H: \Gamma(\lambda)\xi=\lambda^n\xi\}$ 
and $\Gamma(\lambda)=\sum_n\lambda^nE_n$, 
$E_n$ being the projection on $H_n$, $n\in \mathbb Z$.  
In a graded Hilbert module with projections 
$E_n$ as above, one has the commutation relations 
$Z_kE_n=E_{n+1}Z_k$, $1\leq k\leq d$, $n\in\mathbb Z$.  

The graded Hilbert modules that we encounter in this 
paper will all have nonnegative spectrum in the sense 
that $H_n=\{0\}$ for all $n<0$.  
\end{rem}

\begin{rem}[Characterization of Graded Completions]\label{rem0.3}
Graded completions can be characterized abstractly 
as graded Hilbert modules 
$H=H_0\oplus H_1\oplus\cdots$ satisfying 
\begin{enumerate}
\item[A.]$H_0=[v]$ is one-dimensional, 
\item[B.] $Z_1H+\cdots+Z_dH$ is closed, 
\end{enumerate}
for which the map 
$$
f\in\mathcal A\mapsto f(Z_1,\dots,Z_d)v\in H
$$
is injective with dense range.  It follows 
that the homogeneous subspaces 
of such a Hilbert module satisfy 
$$
H_{n+1}=Z_1H_n+\cdots+Z_dH_n, \qquad n=0,1,2,\dots.
$$   
\end{rem}

\section{Standard Hilbert Modules and their Quotients}\label{S:gq}

One can increase the multiplicity of a graded completion 
$G$ to obtain a somewhat more general graded Hilbert module.  
In more detail, let $E$ be a finite-dimensional Hilbert space 
and let $S=G\otimes E$ be the Hilbert module defined by 
$$
f(\xi\otimes\zeta)=(f\cdot\xi)\otimes\zeta, \qquad 
f\in\mathcal A,\quad \xi\in G,\quad \zeta\in E.  
$$
Thus, the coordinate operators $\tilde Z_k$ of $S$ are 
related to the operators $Z_k$ of $G$ by 
$\tilde Z_k=Z_k\otimes\mathbf 1_E$.  
It is convenient 
to ease notation by writing $Z_k$ for $Z_k\otimes\mathbf 1_E$, 
and we usually do so.  We occasionally write $r\cdot G$ for 
$G\otimes\mathbb C^r$.  

\begin{defn}
A {\em standard Hilbert module} is a finite-multiplicity 
version $S=G\otimes E$ of a graded completion $G$.  
\end{defn}

When we want to call attention to the underlying graded completion, 
we say that a standard Hilbert module $S$ is {\em based on $G$} when 
it has the above form $S=G\otimes E$.  Obviously, the direct 
sum of two standard Hilbert modules based on $G$ is a standard Hilbert 
module based on $G$.  On the other hand, it is important to keep in mind that {\em direct sums 
are not allowed across the category of graded completions.}  Indeed, 
if $G_1$ and $G_2$ are two graded completions of $\mathcal A$ that 
are associated with different inner products, then while the direct 
sum $G_1\oplus G_2$ is certainly a graded Hilbert module, it 
need not be a standard Hilbert module (based on any $G$), and 
it may fail to have the favorable properties of standard Hilbert 
modules.  

Standard Hilbert modules carry an obvious grading 
$H=H_0\oplus H_1\oplus\cdots$.  For example, the gauge 
group of $G\otimes E$ is $\Gamma(\lambda)=\Gamma_0(\lambda)\otimes\mathbf 1_E$, 
$\lambda\in\mathbb T$, $\Gamma_0$ being the gauge group of $G$.  
Note too that for any standard Hilbert 
module $S$, the space 
$$
Z_1S+\cdots+Z_dS=(Z_1G+\cdots+Z_dG)\otimes E
$$ 
is closed, and we have 
\begin{equation}\label{gcEq1}
S_{n+1}=Z_1S_n+\cdots+Z_dS_n,\qquad n=0,1,2,\dots.
\end{equation}

Standard Hilbert modules are of degree $0$ in the sense of the 
following general definition.  

\begin{defn}[Degree of a Graded Module]\label{dgDef1}
Let $H=H_0\oplus H_1\oplus H_2\oplus\cdots$ be a graded 
Hilbert module.  The {\em degree} of $H$ is the smallest 
integer $n\geq 0$ such that 
$$
H_{k+1}=Z_1H_k+\cdots+Z_dH_k,\qquad k\geq n.  
$$
If there is no such $n\in\mathbb Z_+$ then the degree 
of $H$ is defined as $\infty$.  
\end{defn}

Hilbert's basis theorem implies that every graded 
submodule of a finitely generated 
graded module over the polynomial algebra $\mathcal A$ has 
a finite number of homogeneous generators.  
It is a straightforward exercise to apply 
that fact to deduce the following:  

\begin{prop}\label{dgProp1}
For every graded submodule $M=M_0\oplus M_1\oplus \cdots $ of a standard Hilbert module 
$S=G\otimes E$ there is an $n=0,1,2,\dots$ such that 
$$
M_{k+1}=Z_1M_k+\cdots+Z_dM_k,\qquad k\geq n.  
$$
\end{prop}

We conclude that {\em Every graded submodule of a standard 
Hilbert module is of finite nonnegative degree.}

\begin{rem}[Degree of Submodules]\label{gqRem2}
A direct application of Proposition \ref{gcProp1} shows   
that the reducing submodules of a standard Hilbert module $G\otimes E$ are 
the submodules $G\otimes F$, where $F$ is a linear subspace 
of $E$.  It follows that a graded submodule $M=M_0\oplus M_1\oplus\cdots$ of 
a standard Hilbert module is a reducing submodule iff
it is of degree $0$.  
More generally, every graded submodule $M$ of $G\otimes E$ 
of degree $n\geq 1$ admits a representation 
$$
M=M_0\oplus M_1\oplus\cdots\oplus M_{n-1}\oplus [\mathcal AM_n], 
$$
where $M_k\subseteq\mathcal A_k\otimes E$  is a linear space of homogeneous 
vector polynomials.   
\end{rem}

Let $S=G\otimes E$ be 
a standard Hilbert module based on a graded 
completion $G$.  We are primarily interested 
in properties of the quotient modules $H=S/M$, where $M$ is a graded submodule of $S$.  
Such quotient modules carry a natural grading $H=H_0\oplus H_1\oplus\cdots$, in 
which $H_k=S_k/M_k$, $k=0,1,2,\dots$.  
A quotient module $H=(G\otimes E)/M$ of this form is called a {\em $G$-quotient}.  
If the submodule $M$ is of degree $n=0,1,2,\dots$, we often 
refer to the quotient as a {\em $G(n)$-quotient}.

\begin{rem}[Ambiguity in the notion of $G(n)$-quotients.]\label{gqRem1} 
Several things must be kept in mind when dealing with 
quotients of standard Hilbert modules.  

1.   Consider first the case in which 
$G$ is the module $H^2$ of the $d$-shift.  The dilation theory of $d$-contractions 
implies that there is an intrinsic characterization of $H^2$-quotients up to 
unitary equivalence: A graded Hilbert module $H=H_0\oplus H_1\oplus\cdots$ is 
unitarily equivalent to an $H^2$-quotient iff its coordinate operators 
$T_1,\dots,T_d$ define a pure finite rank $d$-contraction (we oversimplify 
slightly in order to make the essential point; see \cite{arvSubalgIII} 
and \cite{arvCurv} for more detail).  On the other hand, 
there is no known characterization 
of $H^2$-quotients up to {\em isomorphism}, or even $p$-isomorphism for $p>d$ 
($p$-morphisms are introduced Section \ref{S:pm}).  For 
more general graded completions $G$ in place of $H^2$, very little is known about 
the characterization of $G$-quotients -- even up to unitary equivalence.  

2.  In more explicit terms, given that one knows somehow 
that a graded Hilbert module $H$ is isomorphic 
to a $G$-quotient for some graded completion $G$, one does not know how to 
use the structure of $H$ to obtain information about $G$.  
One does not even know if the structure of some particular $G$-quotient 
$H$ determines $G$ up to isomorphism; though we lack a specific 
example it seems likely that a $G$-quotient 
can be isomorphic to a $G^\prime$-quotient 
when $G$ and $G^\prime$ are graded completions that are not isomorphic.  

3.  Even when one works within the category of $G$-quotients for a {\em fixed} 
graded completion $G$, there is still ambiguity in the notion of $G(n)$-quotient, 
since $n$ is not uniquely determined by the structure of a $G(n)$-quotient.   
Indeed, in the following sections we exploit this ambiguity 
by showing that, up to finite-dimensional  
perturbations, a $G(n)$-quotient for $n\geq 2$ can always be realized up to 
isomorphism as a $G(1)$-quotient - namely a quotient of the form 
$(G\otimes E)/M$ where $M\subseteq G\otimes E$ is a graded submodule 
of {\em degree one}.  
\end{rem}

\section{$p$-essential normality, $p$-morphisms}\label{S:pm}

Let $H$ and $K$ be Hilbert spaces.  
For $p\in[1,\infty]$, we write $\mathcal L^p(H,K)$, 
or more simply $\mathcal L^p$, for the Schatten-von Neumann 
class of all operators $A\in\mathcal B(H,K)$ whose modulus $|T|=\sqrt{T^*T}$ 
satisfies $\tr |T|^p<\infty$ when $p<\infty$, $\mathcal L^\infty$ being 
interpreted as the space of compact operators $\mathcal K(H,K)$.  

\begin{defn}\label{pmDef1}
Let $p\in[1,\infty]$ and let $H$ be a Hilbert module with coordinate 
operators $A_1,\dots,A_d$.  $H$ is said to be $p$-essentially normal 
if the self-commutators $A_jA_k^*-A_k^*A_j$ belong to $\mathcal L^p$ 
for all $1\leq j,k\leq d$.  
\end{defn}

We refer to $\infty$-essentially normal Hilbert modules with 
the shorter term {\em essentially normal}.  
We will make repeated use of the following 
general result which asserts, roughly, that 
submodules and quotients of 
essentially normal Hilbert modules are either very 
good or very bad.  Notice that neither the statement 
nor proof of Proposition \ref{pmProp1} provides information 
about how one might establish the favorable properties.    

\begin{prop}\label{pmProp1}
Let $p\in[1,\infty]$ and let $H$ be a $p$-essentially normal 
Hilbert module 
such that $T_1H+\cdots+T_dH$ is a closed subspace 
of finite codimension in $H$.  For every submodule $M\subseteq H$ 
such that $\{0\}\neq M\neq H$ 
with projection $P=P_M$, 
the following are equivalent:
\begin{enumerate}
\item[(i)] $M$ is $p$-essentially normal.  
\item[(ii)]$H/M$ is $p$-essentially normal.  
\item[(iii)] 
Each commutator $[P,T_1],\dots, [P,T_d]$, belongs to $\mathcal L^{2p}$.  
\item[(iv)] $[P,T_k]^*[P,T_j]\in\mathcal L^p$ for $1\leq j,k\leq d$.  
\item[(v)] $[P,T_j][P,T_k]^*\in\mathcal L^p$ for $1\leq j,k\leq d$.  
\item[(vi)] $M$ is $p$-essentially normal and $T_1M+\cdots+T_dM$ is a 
closed subspace of 
finite codimension in $M$.  
\item[(vii)] $H/M$ satisfies the two conditions of (vi).  
\end{enumerate}
\end{prop}

\begin{proof}
While the proof 
is a straightforward variation of the proof of Proposition 4.1 of \cite{arvpsum}, 
we present the details for completeness.  
Properties (vi) and (vii) obviously imply (i) and (ii), respectively.

(iii)$\iff$(iv)$\iff$(v): Consider the row operator 
whose entries are the commutators 
$C=([P,T_1],[P,T_2],\dots,[P,T_d])$.  $C^*C$ is the 
$d\times d$  matrix $(A_{ij})$ with entries 
$A_{ij}=[P,T_i]^*[P,T_j]$.
After noting that $C^*C\in\mathcal L^p$ iff
$C\in\mathcal L^{2p}$, the equivalence (iii)$\iff$(iv) 
follows.  The proof of (iii)$\iff$(v) is similar.

(i)$\iff$(v):  Letting $B_j$ be the restriction 
of $T_j$ to $M$, we claim 
\begin{equation}\label{psEq2}
[B_j,B_k^*]P=-[P,T_j][P,T_k]^*+P[T_j,T_k^*]P.  
\end{equation}
Indeed, for each $1\leq j,k\leq d$ we can write 
\begin{align*}
[B_j,B_k^*]P=&T_jPT_k^*P-PT_k^*T_jP = 
T_jPT_k^*P-PT_jT_k^*P +P[T_j,T_k^*]P\\
=&-PT_jP^\perp T_k^*P+P[T_j,T_k^*]P.  
\end{align*}
Since $PT_jP^\perp =PT_j-T_jP$, we have 
$PT_jP^\perp T_k^*P=[P,T_j][P,T_k]^*$, 
and (\ref{psEq2}) follows.  Since $H$ is 
$p$-essentially normal we have $[T_j,T_k^*]\in\mathcal L^p$, 
and at this point the equivalence (i)$\iff$(iv) follows from (\ref{psEq2}).  

(i)$\iff$(iv):  Letting $C_j$ be the compression of 
$T_j$ to $M^\perp$, we claim
\begin{equation}\label{psEq3}
[C_j,C_k^*]P^\perp = [P,T_k]^*[P,T_j]+P^\perp [T_j,T_k^*]P^\perp,   
\end{equation}
$P^\perp$ denoting the projection on $M^\perp$.  
Noting that $P^\perp T_j P^\perp = P^\perp T_j$, one has  
\begin{align*}
[C_j,C_k^*]P^\perp 
&=P^\perp T_jT_k^*P^\perp-P^\perp T_k^*P^\perp T_jP^\perp 
\\
&=P^\perp T_k^*T_jP^\perp -P^\perp 
T_k^*P^\perp T_jP^\perp +P^\perp [T_j,T_k^*]P^\perp\\
&=P^\perp T_k^*PT_jP^\perp + P^\perp [T_j,T_k^*]P^\perp.  
\end{align*}
(\ref{psEq3}) follows after one notes that 
$P^\perp T_k^*PT_jP^\perp=[P,T_k]^*[P,T_j]$.  Again, 
we have $[T_j,T_k^*]\in\mathcal L^p$ by hypothesis, 
and the equivalence (i)$\iff$(iv) follows from 
formula (\ref{psEq3}).

(i)$\implies$(iv):  Assuming (i), we have to show that 
$T_1M+\cdots+T_dM$ is closed of finite codimension.  This will follow if we 
show that the restrictions $A_k=T_k\restriction_M$ have the property that 
$A_1A_1^*+\cdots+A_dA_d^*$ is a Fredholm operator.  

For that, we show that this sum has the form $X+K$ where 
$X$ is a positive invertible operator and $K$ is compact.  
Indeed,the hypotheses on $H$ imply 
that $T_1^*T_1+\cdots+T_d^*T_d$ is a positive 
Fredholm operator, which therefore has the form $X_0+F$ 
where $X_0$ is a positive invertible operator and 
$F$ is a finite rank operator.  Hence the compression 
of $X_0+F$ to $M$ has the from $X+K$ 
where $X=P_MX_0\restriction_M$ is a positive invertible 
operator on $M$ and $K$ is the finite rank compression 
of $F$ to $M$.  Using $P_MT_k\restriction_M=T_k\restriction_M$, 
we have  
$$
A_1^*A_1+\cdots +A_d^*A_d=P_M(T_1^*T_1+\cdots+T_d^*T_k)\restriction_M
=P_M(X_0+F)\restriction_M=X+K, 
$$ 
hence $A_1^*A_1+\cdots+A_d^*A_d$ is Fredholm.  Since 
the self-commutators $[A_k^*,A_k]$ are compact by (i), 
it follows that $A_1A_1^*+\cdots+A_dA_d^*$ is a compact 
perturbation of $A_1^*A_1+\cdots+A_d^*A_d$, and is 
therefore Fredholm.

(ii)$\implies$(v): The argument is similar to the one above.  
One notes that the compressions 
$B_k=P_M^\perp T_k\restriction_{M^\perp}$ of the coordinate 
operators to $M^\perp$ 
satisfy 
$$
P_M^\perp(T_1T_1^*+\cdots+T_dT_d^*)\restriction_{M^\perp}=
B_1B_1^*+\cdots+B_dB_d^*,
$$
and uses the fact that $T_1T_1^*+\cdots +T_dT_d^*$ is a positive 
Fredholm operator.  
\end{proof}

Obviously, the notion 
of $p$-essential normality is invariant under unitary equivalence 
of Hilbert modules.  On the other hand, essential normality may 
not be preserved under more general isomorphisms of Hilbert modules - even 
under graded isomorphisms of graded Hilbert modules.  
Some elementary examples are described in Remark \ref{pmRem1} below.  
Similarly, submodules and quotients of essentially normal Hilbert modules need 
not be essentially normal in general (some ungraded examples are exhibited 
in \cite{GRS}).  
The purpose of this section is to establish general conditions 
under which $p$-essential normality does 
propagate to submodules, quotients, and to their images 
under homomorphisms and isomorphisms.  

\begin{defn}[$p$-morphism]\label{pmDef2}
Let $H,K$ be Hilbert modules with respective operator 
$d$-tuples $A_1,\dots,A_d$, $B_1,\dots,B_d$, and 
let $p\in[1,\infty]$.  
By a {$p$-morphism} from $H$ to $K$ we mean 
a homomorphism $L\in\hom(H,K)$ with the additional property  
$L^*B_k-A_kL^*\in\mathcal L^p$ (or equivalently,   
$LA_k^*-B_k^*L\in\mathcal L^p$), $k=1,\dots,d$.  
A {\em $p$-isomorphism} is an isomorphism of Hilbert modules 
that is also a $p$-morphism.  
\end{defn}

One checks readily that the inverse of a $p$-isomorphism 
$L: H\to K$ is a $p$-isomorphism $L^{-1}: K\to H$.  
The following result provides a context in which 
$p$-morphisms of submodules and 
quotient modules arise naturally, and it 
implies that $p$-essential normality of quotient modules 
propagates as desired under appropriate conditions.  

\begin{thm}\label{pmThm1}
Fix $p\in[1,\infty]$.  Let $H$ and $K$ be $p$-essentially normal 
Hilbert modules, let $L: H\to K$ be a $2p$-morphism, and let 
$M\subseteq H$ be a submodule such 
that $L(M)$ is closed.  
If $M$ is $p$-essentially normal then: 
\begin{enumerate}
\item[(i)] $L(M)$ and $M\cap \ker L$ are $p$-essentially normal 
submodules of $K$ and $M$, respectively.
\item[(ii)]The restriction of 
$L$ to $M$ defines a $2p$-morphism 
in $\hom(M,L(M))$ and in $\hom(M,K)$. 
\item[(iii)] The quotients $H/M$ and $K/L(M)$ are 
$p$-essentially normal, and the promoted map  
$\dot L\in \hom(H/M,K/L(M))$ is a $2p$-morphism.    
\end{enumerate}
\end{thm}

\begin{proof}(i): 
Let $X_1,\dots,X_d$ and $Y_1,\dots,Y_d$ be 
the coordinate operators of $H$ and $K$ respectively, 
and consider the operator $B=LP_ML^*\in\mathcal B(K)$.  
Note that the commutators $[B, Y_k]$ must belong 
to $\mathcal L^{2p}$, $k=1,\dots,d$.  Indeed, 
for each $k$ we have $Y_kL=LX_k$, 
$X_kL^*-L^*Y_k\in\mathcal L^{2p}$ because 
$L$ is a $2p$-morphism, 
and Proposition \ref{pmProp1} implies that 
$X_kP_M-P_MX_k\in\mathcal L^{2p}$ because  
$M$ is $p$-essentially normal.  
It follows that 
$$
Y_kB-BY_k=L(X_kP_M-P_MX_k)L^*+LP_M(X_kL^*-L^*Y_k)\in\mathcal L^{2p}.  
$$ 

Note too that, since $N=L(M)$ is closed, 
$B$ is a positive operator with 
closed range $N$, hence either $B$ is invertible 
or $0$ is an isolated point in $\sigma(B)$.  
In either case, there is a positive distance $\epsilon>0$ 
between $0$ and the remaining part $\sigma(B)\cap (0,\infty)$ of 
the spectrum of $B$.  

Let $\Gamma$ be the counter-clockwise oriented curve 
consisting of a rectangle with corners 
$(\epsilon/2,\epsilon/2), (\epsilon/2,-\epsilon/2), 
(\|B\|+\epsilon/2,-\epsilon/2), (\|B\|+\epsilon/2,\epsilon/2)$.  
$\Gamma$ winds once around every point 
of $\sigma(B)\cap (0,\infty)$ and 
has $0$ in its exterior, so the projection $P_N$ can be expressed 
as an operator-valued Riemann integral 
\begin{equation*}
P_N=
\frac{1}{2\pi i}\int_\Gamma R_\lambda \,d\lambda,
\end{equation*}
$R_\lambda$ denoting the resolvent $R_\lambda=(\lambda-B)^{-1}$.  
Perhaps it is appropriate to recall 
that, while in general the right 
side of the displayed formula 
merely defines an idempotent 
with the same range as $B$, in this case 
the integral is a normal 
operator since the $R_\lambda$ are commuting normal 
operators, and a normal idempotent must be 
a self-adjoint projection.  Hence it is $P_N$.

It follows that for any bounded derivation 
of $\mathcal B(K)$ of the form $D(T)=[Y,T]$, 
with $Y\in\mathcal B(K)$, we have 
\begin{equation}\label{pmEq1}
D(P_N)=
\frac{1}{2\pi i}\int_\Gamma D(R_\lambda) \,d\lambda .  
\end{equation}

Notice next that 
\begin{equation}\label{pmEq2}
D(R_\lambda)=R_\lambda D(B)R_\lambda,\qquad \lambda\in \Gamma.  
\end{equation}
Indeed, using $D(S)T+SD(T)=D(ST)$ we have 
$$
D(R_\lambda)(\lambda-B)-R_\lambda D(B)=
D(R_\lambda)(\lambda-B)+R_\lambda D(\lambda-B)=D(\mathbf 1)=0, 
$$
and (\ref{pmEq2}) follows after multiplying on 
the right by $(\lambda-B)^{-1}$.

Now take $Y=Y_k$, $k=1,\dots,d$.  
Using (\ref{pmEq1}) and 
(\ref{pmEq2}), we obtain the formula 
\begin{equation}\label{pmEq3}
[Y_k,P_N]=\frac{1}{2\pi i}\int_\Gamma 
R_\lambda [Y_k,B] R_\lambda\,d\lambda.  
\end{equation}
Since $[Y_k,B]\in\mathcal L^{2p}$ and $\lambda\mapsto R_\lambda$ is 
continuous in the operator norm, it follows that 
$\lambda\in\Gamma\mapsto R_\lambda[Y_k,B]R_\lambda$ 
is a continuous function from $\Gamma$ into 
the Banach space $\mathcal L^{2p}$.  Hence 
(\ref{pmEq3}) expresses $[Y_k,P_N]$ as a Riemann 
integral of a continuous $\mathcal L^{2p}$-valued function, 
and this implies that $[Y_k,P_N]\in\mathcal L^{2p}$.   

Finally, since  
$N$ is a submodule of the $p$-essentially normal 
Hilbert module $K$, 
we may conclude from Proposition \ref{pmProp1} that 
$N$ is a 
$p$-essentially normal submodule of $K$.  

Let $L_0$ be the restriction of $L$ to $M$.  To see 
that $\ker L_0=M\cap \ker L$ is a $p$-essentially 
normal submodule of $H$, consider the 
adjoint 
$L_0^*=P_ML^*$ of $L_0$ as an operator in $B(K,M)$.  
The range of $L_0^*$ is the closed 
subspace $M\ominus (M\cap\ker L)=M\ominus \ker L_0$.  
If we use the 
respective $d$-tuples $Y_1^*,\dots,Y_d^*$ and 
$P_MX_1^*\restriction_M,\dots,P_MX_d^*\restriction_M$ 
to make $K$ and $M$ into $p$-essentially normal 
Hilbert modules over $\mathcal A$, then by 
Proposition \ref{pmProp1}, 
$L_0^*$ becomes a $2p$-morphism in $\hom(K,M)$.  
Thus the argument 
given above implies that the projection on 
$L_0^*(K)$ commutes modulo 
$\mathcal L^{2p}$ with 
$P_MX_1^*\restriction_M,\dots,P_MX_d^*\restriction_M$, 
and hence the projection on $M\cap \ker L$ 
commutes modulo $\mathcal L^{2p}$ with 
$X_1\restriction_M,\dots,X_d\restriction_M$.  
Another application of Proposition 
\ref{pmProp1} shows that 
$M\cap \ker L$ is a $p$-essentially normal Hilbert 
module with respect to the action of 
$X_1\restriction_M,\dots,X_d\restriction_M$.

(ii):  Noting that 
$L_0^*=P_ML^*=P_ML^*P_N$, we have 
\begin{align*}
L_0^*Y_k-X_kL_0^*&=P_ML^*Y_k-X_kP_ML^*
\\&
=P_M(L^*Y_k-X_kL^*)+(P_MX_k-X_kP_M)L^*.  
\end{align*}
The term $L^*Y_k-X_kL^*$ belongs to $\mathcal L^{2p}$ because 
$L$ is a $2p$-morphism, and since $M$ is a $p$-essentially 
normal submodule, Proposition \ref{pmProp1} implies that 
$P_MX_k-X_kP_M\in\mathcal L^{2p}$. 
Hence $L_0$ is a $2p$-morphism in $\hom(M,K)$.  The fact 
that $L_0^*$ is also a $2p$-morphism in $\hom(M,N)$ follows after 
multiplying the previous expressions on the right by 
$P_N$ and arguing the same way.  

(iii): Since $M\subseteq H$ and $N=L(M)\subseteq K$ are 
$p$-essentially normal, Proposition \ref{pmProp1} implies 
that their respective quotients are $p$-essentially normal 
as well.  It remains to show that the promoted map 
$\dot L: H/M\to K/N$ is a $2p$-morphism.   For that, we 
identify $H/M$ with $M^\perp\subseteq H$, $K/N$ with 
$N^\perp\subseteq K$, and $\dot L$ with the map 
$\dot L\xi = P_N^\perp L\xi$, $\xi\in M^\perp$, 
$P_N^\perp$ denoting $\mathbf 1-P_N$.  
We have to 
show that $\dot LX_k^*P_M^\perp-Y_k^*\dot LP_M^\perp$ belongs 
to $\mathcal L^{2p}$, that is,
\begin{equation}\label{pmEq5}
P_N^\perp LX_k^*P_M^\perp-Y_k^*P_N^\perp LP_M^\perp\in\mathcal L^{2p}.  
\end{equation}
The left side of (\ref{pmEq5}) can be written 
$$
P_N^\perp(LX_k^*-Y_k^*L)P_M^\perp + (P_N^\perp Y_k^*-Y_k^*P_N^\perp)LP_M^\perp.  
$$
The term $LX_k^*-Y_k^*L$ belongs to $\mathcal L^{2p}$ because 
$L$ is a $2p$-morphism, and since $N=L(M)$ has been shown to be 
$p$-essentially normal, Proposition \ref{pmProp1} implies that 
$P_N^\perp Y_k^*-Y_k^*P_N^\perp=Y_k^*P_N-P_NY_k^*\in\mathcal L^{2p}$.
(\ref{pmEq5}) follows.   
\end{proof}

\begin{rem}[Estimating the $2p$-norm of $Y_kP_N-P_NY_k$]   
In terms of the operator $L: H\to K$, 
we can estimate the $\mathcal L^{2p}$-norm of 
$[Y_k,P_N]$ from (\ref{pmEq3}) in the obvious way to obtain 
$$
\|[Y_k,P_N]\|_{2p}\leq (2\pi)^{-1}\sup_{\lambda\in\Gamma}\|(\lambda-B)^{-1}\|^2\cdot 
\|[Y_k,B]\|_{2p} \cdot \ell(\Gamma),  
$$
$\ell(\Gamma)$ denoting the length of $\Gamma$.  
Noting that $B=LL^*$, we have 
$\|B\|=\|L\|^2$, hence the length 
of $\Gamma$ is $2\|L\|^2+2\epsilon$;   
and since the minimum distance from $\Gamma$ to $\sigma(Z)$ is $\epsilon/2$,  
we have $\|(\lambda-LL^*)^{-1}\|^2\leq 4/\epsilon^2$,  hence 
\begin{equation}\label{pmEq4}
\|[Y_k,P_N]\|_{2p}\leq \frac{4\|L\|^2+4\epsilon}{\pi\epsilon^2}\|[Z_k,LL^*]\|_{2p}.  
\end{equation}
Here, $\epsilon$ is the length of the gap between 
$0$ and the rest of $\sigma(LL^*)$.  
\end{rem}

\begin{cor}\label{pmCor1}
Fix $p\in[1,\infty]$.  Let $H,K$ be $p$-essentially normal 
Hilbert modules and let $L:H\to K$ be a $2p$-morphism with closed 
range. 

Let $N$ be a (closed) submodule of $L(H)$ and let $M\subseteq H$ 
be its pullback
$$
M=\{\xi\in H: L\xi\in N\}.   
$$
Then $L$ promotes to an isomorphism 
of quotients $\dot L:H/M\to L(H)/N$ 
with the property that 
if either $H/M$ or $L(H)/N$ is $p$-essentially normal, 
then both are $p$-essentially normal.  In that event, 
$\dot L$ is a $2p$-isomorphism.  
\end{cor}

\begin{proof}  In general, 
the promoted map $\dot L\in\hom(H/M,L(H)/N)$ 
is injective because 
$M$ is defined as the full pre-image of $N$ under $L$.  It is 
obviously continuous and surjective; hence the 
closed graph theorem implies that $\dot L$ is an isomorphism 
of Hilbert modules.

Applying Theorem \ref{pmThm1} (i) to $M=H$,  
we see that  $L(H)$ is a $p$-essentially normal submodule of $K$, 
and moreover, $L$ can be viewed as a $p$-morphism in 
$\hom(H,L(H))$.  Thus by replacing $K$ with $L(H)$, we 
can assume that $L$ is surjective.  

If $H/M$ is $p$-essentially normal, then 
$M$ is $p$-essentially normal by Proposition \ref{pmProp1}.  Note that $L(M)=N$.  Indeed, 
the inclusion $\subseteq$ is obvious 
and, since $L$ is surjective, every element $\eta\in N$ has 
the form $\eta=L(\xi)$ for some $\xi\in H$.  Such a 
$\xi$ must belong to $M$ since $M$ is the pre-image of $N$.  
Theorem \ref{pmThm1} (iii) implies that $K/N$ is $p$-essentially normal 
and $\dot L$ is a $2p$-isomorphism.  

Conversely, assume $K/N$ is $p$-essentially normal.  
To prove that $H/M$ is $p$-essentially normal, 
we identify $H/M$ with $M^\perp$, $K/N$ with 
$N^\perp$, and $\dot L: H/M\to K/N$ with the operator 
$$
\dot L: \xi\in M^\perp\mapsto P_N^\perp L\xi\in N^\perp.  
$$
We claim that $L^*(N^\perp)=M^\perp$.  Indeed, 
the preceding paragraph shows that $\dot L$ is an isomorphism 
of $M^\perp$ onto $N^\perp$,  
hence the adjoint $P_M^\perp L^*\restriction_{N^\perp}$ 
of the operator $\dot L$ has range $M^\perp$.    Noting that 
$L^*(N^\perp)\subseteq M^\perp$ simply because 
$L(M)\subseteq N$, it follows that 
$P_M^\perp L^*\restriction_{N^\perp}=L^*\restriction_{N^\perp}$, 
and therefore $L^*(N^\perp)=P_M^\perp L^*(N^\perp)=M^\perp$, 
as asserted.  

The operator $d$-tuples $Y_1^*,\dots,Y_d^*$ and 
$X_1^*,\dots,X_d^*$ make $K$ and $H$ into 
$p$-essentially normal Hilbert modules $K^*$ and $H^*$, 
and $N^\perp$, $M^\perp$ 
become submodules of $K^*$, $H^*$, respectively.  
We can view $L^*$ as a $2p$-morphism of $K^*$ to $H^*$ 
which, by the preceding paragraph, satisfies $L^*(N^\perp)=M^\perp$.   
Now Theorem 
\ref{pmThm1} implies that $M^\perp$ is a $p$-essentially 
normal submodule of $H^*$.  It follows that $H/M$ is 
$p$-essentially normal with respect to its 
original module structure.  
\end{proof}

\begin{rem}[Instability of Essential Normality under Isomorphism]\label{pmRem1}
Essential normality need not be preserved under isomorphism of 
Hilbert modules, even in dimension $d=1$.  As a concrete example, let 
$u_0, u_1,u_2,\dots$ be a sequence of real numbers with the properties
\begin{enumerate}
\item[(i)]$(u_n,u_{n+1})=(0,1)$ for infinitely many  $n=0,1,2,\dots$, 
\item[(ii)]$|u_0+u_1+\cdots+u_n|\leq M<\infty$, $n=0,1,2,\dots$,   
\end{enumerate}
and set 
$$
\lambda_n=e^{u_0+u_1+\cdots+u_n},\qquad n=0,1,2,\dots.     
$$
Let $e_0,e_1,e_2,\dots$ be an orthonormal basis for a Hilbert 
space $H$, let $A\in\mathcal B(H)$ be the simple unilateral shift  
$Ae_n=e_{n+1}$, $n\geq 0$, 
and let $B\in\mathcal B(H)$ be the unilateral weighted shift 
$$
Ae_n=e^{u_{n+1}}e_{n+1},\qquad n=0,1,2,\dots.     
$$
The self-commutator $A^*A-AA^*$ is a rank-one projection, while 
$$
(B^*B-BB^*)e_n=e^{2(u_{n+1}-u_n)}e_n, 
$$ 
so $B$ is not essentially normal by (i).  On the other 
hand, $A$ and $B$ are similar, since one can check  
that the operator $L\in\mathcal B(H)$ defined by 
$Le_n=\lambda_ne_n$ is invertible by (ii), and 
satisfies $LA=BL$.

Note that the operators $A$ and $B$ make 
$H$ into a graded Hilbert module over 
$\mathbb C[z]$ in two ways, and in that case 
$L$ becomes a degree-zero isomorphism 
of graded Hilbert modules that does not preserve 
essential normality.  
\end{rem}

\section{Fredholm Property}\label{S:fp}

 The purpose of this section is to comment on the relation between 
the Fredholm property, essential normality and the cohomology of the 
Koszul complex.  
We also collect an algebraic result for use in Section \ref{S:dg}

\begin{rem}[Koszul complex, Dirac operator]\label{fpRem0}
We briefly recall the definition and most basic properties of 
the Koszul complex and Dirac operator 
of a Hilbert module $H$ of dimension $d$ 
with coordinate operators $T_1,\dots,T_d$.  The 
reader is referred to the original sources \cite{taylor1}, \cite{taylor2} 
for the role of the Koszul complex in operator theory, 
and to   \cite{arvDirac} for more 
on the Dirac operator.  Let $E$ be a Hilbert space of dimension $d$.  For 
each $k=0,1,\dots,d$ there is a Hilbert space of formal $k$-forms 
with coefficients in $H$ 
$$
H\otimes \Lambda^kE
$$
$\Lambda^kE$ denoting the exterior product of $k$ 
copies of $E$ for $k\geq 1$, and $\Lambda^0E$ denoting 
$\mathbb C$.  Fixing an orthonormal basis $e_1,\dots,e_d$ 
for $E$, one obtains creation operators 
$C_1,\dots,C_d: \Lambda^kE\to \Lambda^{k+1}E$, defined uniquely by the requirement 
$$
C_k: \zeta_1\wedge\cdots\wedge \zeta_k\mapsto 
e_k\wedge\zeta_1\wedge\cdots \wedge \zeta_k, \qquad 
\zeta_i\in E,\quad 
k=1,\dots d, 
$$
and a boundary operator $B: H\otimes\Lambda^kE\to H\otimes\Lambda^{k+1}E$ 
by way of 
$$
B = T_1\otimes C_1+\cdots +T_d\otimes C_d.     
$$
Note that $B$ vanishes on $H\otimes \Lambda^dE\cong H\otimes \mathbb C\cong H$.  
This structure gives rise to a complex of Hilbert 
spaces 
\begin{equation}\label{fpEq1}
0\longrightarrow H\longrightarrow H\otimes E
\longrightarrow \cdots\longrightarrow 
H\otimes\Lambda^{d-1}E\longrightarrow H\otimes\Lambda^dE\longrightarrow 0.
\end{equation}
One can assemble the various $H\otimes \Lambda^kE$ into 
a single graded Hilbert space $H\otimes \Lambda E$ by forming 
the direct sum of spaces of $k$-forms 
$$
H\otimes \Lambda E=(H\otimes\Lambda^0E)\oplus\cdots \oplus 
(H\otimes \Lambda^dE),  
$$
thereby making $B$ into a graded operator of degree one 
in $\mathcal B(H\otimes \Lambda E)$ that satisfies  
$B^2=0$.  This structure is called the {\em Koszul complex} 
of $H$.  

The Dirac operator of $H$ is the operator 
$D=B+B^*$.  We are suppressing the Clifford 
structure attached to $D$, which is incidental to our 
needs here.  
The cohomology of the Koszul complex is related to 
the Dirac operator as follows: $B(H\otimes \Lambda^kE)$ has 
finite codimension 
in $\ker (B\restriction_{H\otimes\Lambda^{k+1}E})$, 
for every $k=0,1,\dots,d$
iff $D$ is a self-adjoint Fredholm operator 
\cite{arvDirac}.  
\end{rem}

We say that $H$ has the {\em Fredholm property} 
when these equivalent conditions are satisfied.  
The following general result gives a {\em sufficient} condition 
for Fredholmness that is relatively easy to check in specific 
examples.  While it is a small variation of a result of Curto 
\cite{curFred}, we include a proof for completeness.

\begin{prop}\label{fpProp2}
Let $H$ be an essentially normal 
Hilbert module 
such that $T_1H+\cdots+T_dH$ is a closed subspace 
of finite codimension in $H$. 
Then $H$ has the Fredholm property.  
\end{prop}

\begin{proof}
The hypothesis on the space $T_1H+\cdots+T_dH$ 
is equivalent 
to the assertion that 
$T_1T_1^*+\cdots+T_dT_d^*$ is a 
self-adjoint Fredholm operator.  
We show that the latter property, taken 
together with essential normality, implies that the Dirac operator 
$D$ is Fredholm.  

Let 
$B=T_1\otimes C_1+\cdots+T_d\otimes C_d$ be 
the boundary operator.  Then $D=B+B^*$ and, since 
$T_j$ commutes with $T_k$ and $C_j$ anticommutes with 
$C_k$, one has 
$B^2=0$.  Therefore  
$$
D^2=(B+B^*)^2=B^*B+BB^* = 
\sum_{k,j=1}^d T_k^*T_j\otimes C_k^*C_j+
\sum_{k,j=1}^d T_jT_k^*\otimes C_jC_k^*.  
$$
Using 
$C_jC_k^*=\delta_{jk}\mathbf 1-C_k^*C_j$, we can write 
the second term on the right as 
$$
F\otimes\mathbf 1-\sum_{k,j=1}^d T_jT_k^*\otimes C_k^*C_j, 
$$
where $F=T_1T_1^*+\cdots+T_dT_d^*$, so that 
$$
D^2=F\otimes \mathbf 1 +\sum_{k,j=1}^d [T_k^*,T_j]\otimes C_k^*C_j.  
$$
Since $F\otimes \mathbf 1$ is a Fredholm operator and 
each summand in the second term is compact by hypothesis, it follows 
that $D^2$ is a Fredholm operator.  Since $D$ 
is self-adjoint, $D$ itself 
must be a Fredholm operator.  
\end{proof}

\begin{rem}[Cohomology Type, the Final Three Terms]\label{fpRem0.0}
Let $S=G\otimes\mathbb C^r$ be a standard Hilbert module 
of rank $r$.  If $G$ is essentially normal then so is $S$, 
and therefore Proposition 
\ref{fpProp2} implies that $S$ has the Fredholm property.  
Indeed, it is quite easy to show in that case that the cohomology 
type of $S$ (by which we mean the sequence of Betti numbers 
$(\beta_0,\beta_1,\dots,\beta_d)$ of the Koszul complex of $S$) 
is $(0,0,\dots,0,r)$, the cohomology type 
of the free algebraic module $\mathcal A\otimes\mathbb C^r$ 
of rank $r$.  Thus, {\em essentially normal standard Hilbert modules have the 
same cohomology type as free modules.}

In particular, the behavior of the boundary operator 
of the Koszul complex of a standard Hilbert module is 
specified at every stage $S\otimes\Lambda^kE$, 
$0\leq k\leq d$.  In 
this paper we shall only have to 
refer to the last three terms of the Koszul complex 
\begin{align*}
B_{d-2}&:S\otimes \Lambda^{d-2}E\to S\otimes \Lambda^{d-1}E,\\
B_{d-1}&: S\otimes \Lambda^{d-1}E\to S\otimes \Lambda^dE=S,   
\end{align*}
for graded Hilbert modules whose last two 
Betti numbers are $0,r$.  In such cases, 
the two boundary operators satisfy a) 
$\ran B_{d-2}=\ker B_{d-1}$ and b) $\ran B_{d-1}$ {\em is a closed 
subspace of $S$ of finite codimension $r$}.  

We now point out how 
the assertions a) and b) can be cast into a more concrete form 
for graded modules $H=H_0\oplus H_1\oplus\cdots$ 
with coordinate operators $T_1,\dots,T_d$,  
that satisfy  $0<\dim H_0<\infty$ and 
$$
H_{n+1}=T_1H_n+\cdots+T_dH_n, \qquad n=0,1,2,\dots.   
$$

Let us first consider b).  Note that $H\otimes\Lambda^{d-1}E$ 
can be identified 
with the direct sum $d\cdot H$ of $d$ copies of 
$H$ in such a way that, up to sign, the boundary map 
$B_{d-1}:H\otimes\Lambda^{d-1}E\to H$ becomes identified 
with the row operator 
$$
\bar T: (\xi_1,\dots,\xi_d)\in d\cdot H\mapsto 
T_1\xi_1+\cdots+T_d\xi_d\in H.  
$$
The range 
of this map 
is dense in $H_1\oplus H_2\oplus\cdots = H_0^\perp$ 
and will be of finite (algebraic) 
codimension in $H$ $\iff$ $T_1H+\cdots+T_dH$ is 
closed $\iff$ $T_1T_1^*+\cdots +T_dT_d^*$ 
is a self-adjoint Fredholm operator.  
Note too that 
since 
the orthocomplement of 
$\bar T(d\cdot H)$ is the space 
$H_0\neq \{0\}$, it follows that 
the complex (\ref{fpEq1}) cannot 
be exact at the final term $H\otimes \Lambda^dE$.

In more explicit terms, a) 
makes the assertion that 
for every $d$-tuple of vectors $\xi_1,\dots,\xi_d$ in $H$ 
with the property $T_1\xi_1+\cdots+T_d\xi_d=0$, there is a skew-symmetric 
array $\eta_{ij}=-\eta_{ji}\in H$, $1\leq i,j\leq d$, such that 
\begin{equation}\label{fpEq3}
\xi_k=\sum_{j=1}^dT_j\eta_{jk},\qquad k=1,2,\dots,d.  
\end{equation}

We require the representation 
(\ref{fpEq3}) in Proposition \ref{dgProp3} below.   For now, we simply note that 
if $\xi_1,\dots,\xi_d\in H_n$ are homogeneous 
vectors of degree 
$n=0,1,2,\dots$ satisfying $T_1\xi_1+\cdots+T_d\xi_d=0$, 
then since each operator 
$T_j$ is of degree one, 
the vectors $\eta_{ij}$  of (\ref{fpEq3}) 
must belong to $H_{n-1}$.  
In particular, for the case $n=0$ we have $H_{-1}=\{0\}$, and therefore 
\begin{equation}\label{fpEq2}
\xi_1,\dots,\xi_d\in H_0, \quad T_1\xi_1+\cdots +T_d\xi_d=0 
\implies \xi_1=\cdots=\xi_d=0.  
\end{equation}
\end{rem}

\section{Kernels of Degree $1$}\label{S:dg}

Let $S=G\otimes E$ be a standard Hilbert module and let 
$L:d\cdot S\to S$ be the row operator 
\begin{equation}\label{dgEq1}
L(\xi_1,\dots,\xi_d)=Z_1\xi_1+\cdots+Z_d\xi_d, 
\qquad \xi_1,\dots,\xi_d\in S.  
\end{equation}
We have seen that the range of $L$ is a closed submodule 
of $S$, $\ran L=S_1\oplus S_2\oplus\cdots$, and in particular 
$\ran L$ is of degree $1$.  The kernel of $L$ 
is a graded submodule of $d\cdot S$ that will occupy a 
central position in the following section.  In this section 
we calculate its degree.

\begin{prop}\label{dgProp3}
The kernel 
$K=\ker L$ of the row operator of (\ref{dgEq1}) 
is a graded submodule of $d\cdot S$ of degree $1$.  
\end{prop}

\begin{proof}
That $K=K_0\oplus K_1\oplus\cdots$ is graded is clear from the 
fact that $L$ maps $(d\cdot S)_n$ into $S_{n+1}$ for every $n=0,1,2,\dots$, 
and (\ref{fpEq2}) asserts that $K_0=\{0\}$.  
Hence $K=K_1\oplus K_2\oplus\cdots$.  

Since the inclusion $Z_jK_n\subseteq K_{n+1}$, $1\leq j\leq d$, is obvious,   
we have to prove 
\begin{equation*}
K_{n+1}\subseteq Z_1K_n+\cdots+Z_dK_n,\qquad n\geq 1.  
\end{equation*}

To that end, choose $\xi\in K_{n+1}$, say 
$\xi=(\xi_1,\dots,\xi_d)$ with $\xi_j\in S_{n+1}$ 
satisfying 
$Z_1\xi_1+\cdots+Z_d\xi_d=0$.  Formula (\ref{fpEq3}) implies 
that 
there is an antisymmetric array 
$\eta_{ij}=-\eta_{ji}\in S_n$, $1\leq i,j\leq d$, 
such that $\xi_j=\sum_iZ_i\eta_{ij}$ for $j=1,\dots,d$.  
We consider the array  $(\eta_{ij})_{1\leq i<j\leq d}$ as 
an element of $q\cdot S$, where $q=d(d-1)/2$.  This 
$q$-tuple of vector polynomials 
is homogeneous of degree $n$, 
and therefore it belongs 
to $Z_1(q\cdot S)_{n-1}+\cdots+Z_d(q\cdot S)_{n-1}$.  
Thus we can find a set of $d$ such arrays 
$(\zeta_{ij}^1)_{1\leq i<j\leq d},\dots ,(\zeta_{ij}^d)_{1\leq i<j\leq d}\in q\cdot S$, 
each component of which is homogeneous of degree $n-1$, 
such that 
$$
\eta_{ij}=\sum_{p=1}^d Z_p\zeta_{ij}^p,\qquad 1\leq i<j\leq d.  
$$ 
Setting $\zeta_{ij}^p=-\zeta_{ji}^p$ for $i>j$ and 
$\zeta_{ii}^p=0$ for $1\leq i\leq d$, we obtain antisymmetric 
arrays $(\zeta_{ij}^1),\dots,(\zeta_{ij}^d)$ which, since 
$\eta_{ij}$ is also antisymmetric, satisfy 
$$
\eta_{ij}=\sum_{p=1}^dZ_p\zeta_{ij}^p,\qquad 1\leq i,j\leq d.  
$$
It follows that 
$$
\xi_i=\sum_{j,p=1}^dZ_jZ_p\zeta_{ij}^p, 
$$
and therefore $\xi$ can be written as a linear combination 
\begin{equation}\label{dgEq2}
\xi=(\xi_1,\dots,\xi_d)=Z_1\omega_1+\cdots +Z_d\omega_d
\end{equation}
of elements 
$$
\omega_p=(\sum_j Z_j\zeta_{1j}^p, 
\cdots, \sum_jZ_j\zeta_{dj}^p), \qquad p=1,\dots,d.  
$$
Each $\omega_p$ belongs to $d\cdot S_n$, 
and we have 
$$
L\omega_p=\sum_{i,j=1}^d Z_iZ_j\zeta_{ij}^p=0
$$
because for $p$ fixed, 
$Z_iZ_j\zeta_{ij}^p$ is antisymmetric in $i$ and $j$.  
Hence $\omega_p\in K_n$, 
so that (\ref{dgEq2}) 
exhibits $\xi$ as an element of $Z_1K_n+\cdots+Z_dK_n$.  
\end{proof}

\section{Stability of Quotients under Shifting}\label{S:st}

One can shift a graded Hilbert module $H=H_0\oplus H_1\oplus\cdots$ 
to the left to obtain another graded Hilbert module 
$H^\prime=H_1\oplus H_2\oplus \cdots$, in which the grading on $H^\prime$ 
is given by $H^\prime_n=H_{n+1}$, $n\geq 0$.  One can also 
view $H^\prime$ as submodule of $H$ of
codimension $\dim H_0$, with a grading different 
from the inherited grading.  
In this section we show 
that in an appropriate sense, 
the class of graded Hilbert modules 
that are isomorphic to quotients $(G\otimes E)/M$ of standard Hilbert 
modules based on a fixed graded completion $G$ is stable 
under this operation of left-shifting.  

Example \ref{} shows that in general, isomorphisms of 
Hilbert modules may not preserve essential normality.  
We first strengthen the notion of isomorphism so as to 
get rid of this anomaly.  

\begin{defn}[Strong Equivalence]\label{stDef1}
Let $H$ and $K$ be Hilbert modules with respective 
operator $d$-tuples $X_1,\dots,X_d$ and 
$Y_1,\dots, Y_d$.  $H$ and $K$ are 
said to be {\em strongly equivalent} 
if there is an isomorphism of 
Hilbert modules $L:H\to K$ such that 
$LX_k^*-Y_k^*L$ is compact for every $k=1,\dots,d$.  

In the category of graded Hilbert modules, 
strong isomorphisms are required to preserve 
degree: $L(H_n)=K_n$, $n=0,1,\dots$.  
\end{defn}

\begin{rem}[Strong Equivalence and Essential Normality]
\label{stRem1}
Strong equivalence is obviously an equivalence 
relation.  More significantly, if an essentially normal 
Hilbert module $H$ is strongly equivalent to a Hilbert 
module $K$, then $K$ must be essentially normal as well.  
Indeed, 
choosing $L:H\to K$ as in Definition \ref{stDef1}, 
it follows that $[X_k,L^*L]=X_kL^*L-L^*LX_k$ is compact, hence 
$[X_k,(L^*L)^{-1}]$ is compact, and writing 
$Y_k=LX_kL^{-1}$, we find that  for $1\leq j,k\leq d$, 
\begin{align*}
Y_kY_j^*-Y_j^*Y_k&=LX_k(L^*L)^{-1}X_j^*L^*-L^{*-1}X_j^*L^*LX_kL^{-1}
\\&=L[X_k,(L^*L)^{-1}]+L^{*-1}(X_kX_j^*-X_j^*X_k)L^{-1}\in\mathcal K.  
\end{align*}
\end{rem}

\begin{thm}[Stability of $G$-quotients]
\label{stThm1} Let $G$ be an essentially normal graded completion 
and let 
$H=H_0\oplus H_1\oplus\cdots$ be a graded Hilbert module 
that is strongly 
equivalent to a $G(n)$-quotient for some $n\geq 2$.  Then the shifted 
module $H^\prime=H_1\oplus H_2\oplus\cdots$ is strongly equivalent
to a $G(n-1)$-quotient.  
\end{thm}

\begin{proof}
By hypothesis, there is a standard Hilbert module 
$S=G\otimes E$ based on $G$ and a graded submodule 
$M\subseteq S$ of degree $n$ such that $H$ is 
strongly equivalent to $S/M$.  By Remark \ref{stRem1} 
we may, without loss of 
generality, take $H=S/M$ and $H_k=S_k/M_k$, 
$k=0,1,2,\dots$.  Letting $Z_1,\dots,Z_d$ be 
the coordinate operators of $S$, 
one sees that the shifted module $H^\prime$ is 
given by 
$$
H^\prime=S_1/M_1\oplus S_2/M_2\oplus\cdots=
(Z_1S+\cdots+Z_dS)/(M_1\oplus M_2\oplus\cdots).  
$$
We have to show that there is another standard 
Hilbert module $S^\prime=G\otimes E^\prime$ based 
on $G$, 
a graded submodule $M^\prime\subseteq S^\prime$ 
of degree $n-1$, and a strong isomorphism 
$L$ that maps $S^\prime/M^\prime$ to $H^\prime$.  

Let $S^\prime=d\cdot S= G\otimes d\cdot E$ be the direct sum of 
$d$ copies of $S$, and let $L:d\cdot S\to S$ be the 
row operator $(Z_1,\dots,Z_d)$:
$$
L(\xi_1,\dots,\xi_d)=Z_1\xi_1+\cdots+Z_d\xi_d,\qquad 
\xi_1,\dots,\xi_d\in S.
$$  
Since $G$ is a graded completion, the range of $L$ 
is the closed submodule of $S$ given 
by $Z_1S+\cdots +Z_dS=S_1\oplus S_2\oplus\cdots$.  
Let $M^\prime$ be the pullback of $M$
$$
M^\prime=\{\zeta\in d\cdot S: L\zeta\in M\}.  
$$
We have $L(S^\prime_k)=S_{k+1}$ for $k=0,1,2,\dots$, 
hence $M^\prime=M^\prime_0\oplus M^\prime_1\oplus\cdots$ 
is a graded submodule 
of $S^\prime=d\cdot S$ that contains the kernel of 
$L$ and satisfies 
\begin{equation}\label{stEq3}
L(M^\prime_k)=M_{k+1},\qquad k=0,1,2,\dots.   
\end{equation}
Therefore  
$L(M^\prime)=M\cap L(S^\prime)=M_1\oplus M_2\oplus\cdots$.   
It follows that $L$ promotes to an 
isomorphism of Hilbert modules 
$\dot L: S^\prime/M^\prime \to H^\prime$.  

Since $G$ is essentially normal 
the commutators $Z_j^*Z_k-Z_kZ_j^*$ are compact; 
therefore since  
$L=(Z_1,\dots,Z_d)$ is a row 
operator whose components $Z_k$ commute modulo $\mathcal K$ 
with $Z_1^*,\dots,Z_d^*$, $L$ 
must be an 
$\infty$-morphism.  At this point we can apply 
the last sentence of 
Corollary \ref{pmCor1} to conclude that $\dot L$ 
implements a strong isomorphism of graded Hilbert modules
$S^\prime/M^\prime\cong H^\prime$.  

It remains to show that the pullback 
$M^\prime$ is of degree $n-1$, i.e., 
\begin{equation}\label{stEq1}
Z_1M^\prime_k+\cdots+Z_dM^\prime_k=M^\prime_{k+1}, 
\qquad k\geq n-1,
\end{equation}
and 
\begin{equation}\label{stEq2}
Z_1M^\prime_{n-2}+\cdots+Z_dM^\prime_{n-2}\neq M^\prime_{n-1}.  
\end{equation}

The inclusion $\subseteq$ of (\ref{stEq1}) follows 
from the fact that $M^\prime$ is a graded submodule
of $S^\prime$.  
To prove the opposite inclusion, 
let $k\geq n-1$ and 
choose $\xi\in M^\prime_{k+1}$.  Using (\ref{stEq3}) 
and the fact that $M$ is of degree $n$, we have 
\begin{align*}
L\xi\in M_{k+2}&=Z_1M_{k+1}+\cdots+Z_dM_{k+1}
=Z_1L(M^\prime_{k})+\cdots+Z_dL(M^\prime_k)
\\
&=L (Z_1 M^\prime_{k}+\cdots+Z_d M^\prime_k), 
\end{align*}
and therefore $\xi\in Z_1 M^\prime_{k}+\cdots+Z_d M^\prime_k+\ker L$.  
Noting 
that $\ker L$ is a graded submodule  
$$
\ker L=K_1\oplus K_2\oplus \cdots, 
$$ 
 and that both 
$\xi$ and the space $Z_1 M^\prime_k+\cdots+Z_dM^\prime_k$ 
are homogeneous of 
degree $k+1$, we must have $\xi\in (Z_1M^\prime_k+\cdots+Z_dM^\prime_k)+K_{k+1}$.  
Since $k+1\geq n\geq 2$, 
Proposition \ref{dgProp3} implies 
that $K_{k+1}=Z_1 K_{k}+\cdots+Z_dK_k$, 
so that 
$$
\xi\in Z_1(M^\prime_k+K_k)+\cdots+Z_d(M^\prime_k+K_k).
$$  
Finally, since $M^\prime$ contains 
$\ker L$ we have $K_k\subseteq M^\prime_k$, 
so that 
$$
\xi\in Z_1(M^\prime_k+K_k)+\cdots+Z_d(M^\prime_k+K_k)
=Z_1M^\prime_k+\cdots+Z_dM^\prime_k, 
$$ 
and (\ref{stEq1}) is established.  

To prove (\ref{stEq2}) we proceed contrapositively.  If (\ref{stEq2}) 
fails then 
$$
M^\prime_{n-1}=Z_1M^\prime_{n-2}+\cdots+Z_dM^\prime_{n-2}
$$
and, after applying $L$ and using $LZ_j=Z_jL$ for 
$j=1,\dots,d$, we obtain
\begin{align*}
M_{n-1}&=L(M^\prime_{n-1})=Z_1L(M^\prime_{n-2})+\cdots+
Z_dL(M^\prime_{n-2})\\
&=Z_1M_{n-2}+\cdots+Z_dM_{n-2}.  
\end{align*}
This implies that $\deg M\leq n-1$, contradicting 
$\deg M=n$.  
\end{proof}

\section{Linearizations}\label{S:li}

We now assemble the preceding 
results to deduce the linearization result alluded to in 
the introduction and discuss its implications for the issue 
of essential normality.    Throughout the section, 
$G$ will denote an essentially normal graded completion of 
$\mathcal A$.  

\begin{thm}[Linearization of $G$-quotients]\label{liThm1}
Let $M\subseteq G\otimes E$ be a graded submodule of a 
standard Hilbert module based on $G$, 
let $H=(G\otimes E)/M$ be its 
quotient, and consider the natural 
grading of $H=H_0\oplus H_1\oplus\cdots$.  
If the degree $n=\deg M$ of $M$ is larger than $1$, 
then the shifted module 
\begin{equation}\label{liEq1}
H_{n-1}\oplus H_{n}\oplus H_{n+1}\oplus\cdots 
\end{equation}
is strongly equivalent to the quotient $S^\prime/M^\prime$ 
of another standard Hilbert module $S^\prime=G\otimes E^\prime$ 
based on $G$ by a graded submodule $M^\prime$ of degree $1$.  
\end{thm}

\begin{proof}
The Hilbert module (\ref{liEq1}) is seen to be strongly equivalent 
to a $G(1)$-quotient after a straightforward iteration 
of Theorem \ref{stThm1} through $n-1$ steps.  We omit 
the details.  
\end{proof}

It is convenient to introduce the following terminology.  

\begin{defn}\label{liDef1}
Let $H=(G\otimes E)/M$ be a $G$-quotient.  
By a {\em linearization} of $H$ we mean 
a $G$-quotient $S^\prime/M^\prime$ with the following properties:
\begin{enumerate}
\item[(i)] $M^\prime$ is a degree $1$ graded submodule of $S^\prime$.  
\item[(ii)]$S^\prime/M^\prime$ is strongly equivalent to a closed submodule 
of finite codimension in $H$.  
\end{enumerate}
\end{defn}

\begin{rem}[Nonuniqueness of Linearizations]\label{liRem1}
Let $H=H_0\oplus H_1\oplus\cdots$ be a $G(n)$ quotient for 
some $n\geq 2$.  Theorem \ref{liThm1} asserts that the 
$n-1$-shifted submodule (\ref{liEq1}) is strongly equivalent to 
a $G(1)$-quotient $(G\otimes E^\prime)/M^\prime$.  If 
some other $G(1)$-quotient $(G\otimes E^{\prime\prime})/M^{\prime\prime}$ 
is strongly equivalent to the module of (\ref{liEq1}), then of course 
the two $G(1)$-quotients $(G\otimes E^\prime)/M^\prime$ 
and $(G\otimes E^{\prime\prime})/M^{\prime\prime}$ must be 
strongly equivalent.  Thus one would like to know 
the extent to which a $G(1)$-quotient is uniquely determined 
by its strong equivalence class.  In the case where 
$G$ is the space $H^2$ of the $d$-shift, one knows 
that the unitary equivalence class of quotient modules 
such as $(H^2\otimes E)/M$ uniquely determines 
the pair $M\subseteq H^2\otimes E$ up to unitary 
equivalence once one 
imposes a natural minimality requirement  
\cite{arvSubalgIII}.  On the other hand, even for this 
case $G=H^2$, we no of no classification of 
$G$-quotients (even $G(1)$-quotients) up to 
strong equivalence.  

In a sense, Theorem \ref{liThm1} provides the 
simplest linearization of a $G(n)$-quotient for $n\geq 2$.  
Indeed, one can ask what happens 
when one shifts the module of (\ref{liEq1}) further 
to the left.  It is not hard to show that further 
shifting does {\em not} reduce $G(1)$-quotients to $G(0)$-quotients.  
Rather, it simply leads to an infinite sequence of 
new Hilbert 
modules 
$$
H_k\oplus H_{k+1}\oplus H_{k+2}\oplus\cdots, \qquad 
k=n-1,n, n+1,\dots,  
$$
each of which is strongly equivalent to a $G(1)$-quotient.  
We omit the proof since we do not require that result.  
However, from these remarks it is apparent that {\em every $G$-quotient 
has infinitely many linearizations}.  
\end{rem}

Theorem \ref{liThm1} 
has the following consequence, 
which reduces the problem of establishing essential 
normality of quotients to that of establishing essential 
normality of linearized quotients:  

\begin{cor}\label{liCor1}
Let $G$ be an essentially 
normal graded completion, and let $H=(G\otimes E)/M$ be a 
$G$-quotient.  Then $H$ has linearizations, and the 
following are equivalent:
\begin{enumerate}
\item[(i)]$H$ is essentially normal.  
\item[(ii)] Some linearization of $H$ is essentially normal.  
\item[(iii)] Every linearization of $H$ is essentially normal.  
\end{enumerate} 
\end{cor}

\begin{proof}  Let $M\subseteq H$ be a closed submodule 
of a Hilbert module such that $H/M$ is finite-dimensional.  
Since the projection on $M$ is a finite-rank perturbation 
of the identity, Proposition \ref{pmProp1} implies that $H$ is essentially 
normal iff $M$ is essentially normal.  Thus, the existence 
of linearizations follows from the statement of 
Theorem \ref{liThm1}.  
Since strong equivalence preserves essential normality 
of Hilbert modules, the equivalence of (i) -- (iii) follows.  
\end{proof}

\section{Structure of Linearized Quotients}\label{S:lq}

Linearized quotients admit a particularly concrete 
description rooted in basic algebra, especially 
in the case where $G$ has maximum symmetry in the sense 
of Appendix A.  In particular, the problem of establishing 
essential normality for graded quotients $(H^2\otimes E)/M$ 
is reduced to the problem of establishing essential normality 
for these concrete examples.  We spell that out in this section, 
concentrating on the maximally symmetric cases.  
It is exasperating that the problem of 
establishing essential normality for these 
examples remains open.

Throughout the section, $S=G\otimes E$ 
denotes a standard Hilbert module based 
on a graded completion $G$ of the polynomial algebra $\mathcal A$.  
Note first that $\mathcal A\otimes E$ is a graded module over the 
polynomial algebra $\mathcal A$ that has $S$ as it closure.  
More explicitly, the tensor product 
of the $G$-inner product on $\mathcal A$ and the given 
inner product 
on $E$ gives rise to an inner product on $\mathcal A\otimes E$, 
and the completion of $\mathcal A\otimes E$ in that norm 
is the Hilbert module $S=G\otimes E$.  
Throughout the section, it will be convenient 
to view elements $f$ of 
$\mathcal A\otimes E$ as $E$-valued 
vector polynomials 
$$
z\in\mathbb C^d\mapsto f(z)\in E .
$$
The homogeneous summand $S_n$ is the 
space $\mathcal A_n\otimes E$ of homogeneous polynomials 
of degree $n$ that take values in $E$, $n=0,1,2,\dots$.  

For every linear subspace $V\subseteq E$, 
the closure of 
$\mathcal A\otimes V$ in $S$ is a (necessarily graded) reducing 
submodule of $S$; and by Remark \ref{gqRem2}, every 
reducing submodule of $S$ is obtained in this way 
from a uniquely determined subspace $V\subseteq E$.   We begin 
by summarizing the elementary features 
of the ``pointwise" description of these 
reducing submodules.

\begin{prop}\label{lqProp1}
Let $V$ be a linear subspace of $E$.  The space $\mathcal V$ of all 
vector polynomials $f\in \mathcal A\otimes E$ that satisfy
\begin{equation}\label{lqEq0}
f(z)\in V,
{\ \rm for\ all\ } z=(z_1,\dots,z_d)\in \mathbb C^d  
\end{equation}
is graded in that 
$\mathcal V=\mathcal V_0\dotplus \mathcal V_1\dotplus\mathcal V_2\dotplus\cdots$, 
where $\mathcal V_n=\mathcal A_n\otimes V$ consists of 
all homogeneous polynomials of degree $n$ that satisfy (\ref{lqEq0}).  
The closure of $\mathcal V$ in $S=G\otimes E$ is the reducing 
submodule $G\otimes V$.  
\end{prop}

\begin{proof}  Let $f\in\mathcal A\otimes E$ be a 
polynomial satisfying (\ref{lqEq0}).  Then for each $r\geq 0$, the 
scaled polynomial $f_r(z)=f(rz)$ also satisfies 
(\ref{lqEq0}).  
After noting the Taylor expansion of $f$
$$
f(rz)=f_0(z)+rf_1(z)+r^2f_2(z)+\cdots,\qquad z\in\mathbb C^d,
$$
and carrying out obvious differentiations 
with respect to $r$, one 
finds that $f_n(z)\in V$, $z\in \mathbb C^d$, $n\geq 0$,
hence each homogeneous polynomial $f_n$ 
satisfies (\ref{lqEq0}).  
By the preceding remarks, 
homogeneous polynomials of degree $n$ that satisfy (\ref{lqEq0}) 
are the elements of $\mathcal A_n\otimes V$.  
Since $f_0+f_1+\cdots = f$, it follows that $f$ belongs 
to $\mathcal V_0+\mathcal V_1+\cdots$.  

It follows that the closure of $\mathcal V$ in the norm 
of $S=G\otimes E$ is a graded submodule with homogeneous 
summands $\mathcal V_n=\mathcal A_n\otimes V$, $n=0,1,2,\dots$.  Hence 
it is $G\otimes V$.  
\end{proof}

\begin{rem}[The differential operator 
$D_G:\mathcal A\otimes E\to \mathcal A\otimes d\cdot E$]\label{lqRem1}
Fix a graded completion $G$ and 
consider the adjoints $Z_k^*$, $1\leq k\leq d$, of the coordinate operators 
of $G$.  Recalling that $G_n=\mathcal A_n$ and that $Z_k^*G_n\subseteq G_{n-1}$ 
for $n\geq 1$ and $Z_k^*G_0=\{0\}$, it follows that each $Z_k^*$ leaves 
$\mathcal A$ invariant, and carries $\mathcal A_n$ into $\mathcal A_{n-1}$.   
It goes without saying that, while the action of $Z_k$ on $\mathcal A$ 
depends only on the algebraic structure of $\mathcal A$ (and the specified 
basis for $\mathcal A_1$), the 
operators $Z_k^*$ depend strongly on the inner product 
associated with $G$.  

For higher multiplicity standard Hilbert modules 
$S=G\otimes E$ based on $G$, the adjoints of 
the coordinate operators of $S$ are the higher multiplicity versions  
$Z_k^*\otimes \mathbf 1_E$ of the corresponding operators on $G$.    
We will continue to abuse notation 
by writing $Z_k$ and $Z_k^*$ for the coordinate operators of $S$ 
and their adjoints, whenever it does not lead to confusion.  Thus, 
{\em each $Z_k^*$ acts as a linear operator of degree $-1$ on the 
graded algebraic module $\mathcal A\otimes E$}.  We find it convenient 
to think of $Z_k^*$ as a generalized first order differential operator  
(for example, see (\ref{lqEq0.2}) below).    

If we realize $d\cdot(\mathcal A\otimes E)=\mathcal A \otimes d\cdot E$ 
as the space of $d$-tuples $(f,\dots,f_d)$ of vector polynomials 
$f_k\in\mathcal A\otimes E$, 
then we can define a linear operator 
$D_G:\mathcal A\otimes E\to \mathcal A\otimes d\cdot E$
as follows:
\begin{equation}\label{lqEq0.1}
D_Gf=(Z_1^*f,\dots,Z_d^*f),\qquad f\in\mathcal A\otimes E.  
\end{equation}
The case where $G=H^2$ is the space of the $d$-shift 
is noteworthy, since in that case the adjoints of the coordinate 
operators of $H^2$ act on polynomials by way of 
\begin{equation}\label{lqEq0.2}
Z_k^*f(z)=(N+\mathbf 1)^{-1}\frac{\partial f}{\partial z_k}(z), 
\qquad z\in\mathbb C^d.  
\end{equation}
where $N$ is the number operator of $H^2$.  It follows that the 
action of $D_{H^2}$ on vector polynomials in $\mathcal A\otimes E$ 
is given by 
\begin{equation}\label{lqEq1}
D_{H^2}f(z)=(N+\mathbf 1)^{-1}
(\frac{\partial f}{\partial z_1}(z)\dots,\frac{\partial f}{\partial z_d}(z))
=(N+\mathbf 1)^{-1}\nabla f(z) 
\end{equation}
$\nabla$ denoting the classical gradient operator.  
More generally, if $G$ is a graded completion associated with a maximally  
symmetric inner product as in Appendix A, then $D_G$ takes the 
form 
\begin{equation}\label{lqEq1.1}
D_Gf(z)=u(N)\nabla f(z),\qquad z\in\mathbb C^d,
\end{equation}
where $u$ is an appropriate bounded function in $C[0,\infty)$.  
\end{rem}

In order to keep the statement of results 
as simple as possible, we confine attention 
to maximally symmetric graded completions $G$, so that the operators 
$D_G: \mathcal A\otimes E\to \mathcal A\otimes d\cdot E$ 
have the form (\ref{lqEq1.1}).  
We seek a concrete description of all $G(1)$-quotients 
in elementary terms  (Theorem \ref{lqThm1}).  At the end of the section 
we indicate how that description should be modified for more 
general graded completions.  

We begin by giving a purely algebraic description of a family 
of subspaces of $\mathcal A\otimes E$ that are invariant 
under the action of the differentiation operators 
$\frac{\partial}{\partial z_1},\dots\frac{\partial}{\partial z_d}$.  
Let $E$ be a finite-dimensional vector space, let 
$V$ be a linear subspace 
of the direct sum $d\cdot E$ of $d$ copies of $E$,  and consider the 
space $\mathcal E_V$ of all polynomials $f\in\mathcal A\otimes E$ with 
the property:
\begin{equation}\label{lqEq2}
\nabla f(z)\in V,\qquad z=(z_1,\dots,z_d)\in \mathbb C^d.  
\end{equation}
One finds that 
$\mathcal E_V$ is a graded vector space 
$$
\mathcal E_V=\mathcal E_V(0)\dotplus\mathcal E_V(1)\dotplus\mathcal E_V(2)\dotplus\cdots, 
$$
in which $\mathcal E_V(n)$ denotes the space of all homogeneous polynomials 
of degree $n$ that satisfy (\ref{lqEq2}) -- 
the argument being a minor variation on the proof of Proposition \ref{lqProp1}.  
One also has 
\begin{equation}\label{lqEq6}
\frac{\partial}{\partial z_k}\mathcal E_V(n)\subseteq \mathcal E_V(n-1),\qquad n=1,2,\dots, 
\quad k=1,\dots, d, 
\end{equation}
and $\mathcal E_0$, the space of constant polynomials, is annihilated 
by $\frac{\partial}{\partial z_1}, \dots,\frac{\partial}{\partial z_d}$.  

\begin{thm}\label{lqThm1}
Let $S=G\otimes E$ be a standard Hilbert module based on 
a maximally symmetric graded completion of $\mathcal A$ 
and let $Z_1^*,\dots,Z_d^*$ be the adjoints of the 
coordinate multiplications of $S$.  
Then $Z_k^*\mathcal E_V\subseteq \mathcal E_V$, 
hence the orthocomplement of $\mathcal E_V$ in $S$ 
is a graded submodule of $S$, and the 
closure of $\mathcal E_V$ in $S$ is identified with the quotient module 
\begin{equation}\label{lqEq3}
H_V=S/\mathcal E_V^\perp .  
\end{equation}
This quotient has 
the property that its submodule $M=\mathcal {E}_V^\perp$ 
is a degree $1$ graded submodule of $Z_1S+\cdots+Z_dS$.  

Conversely, every degree $1$ graded submodule $M\subseteq Z_1S+\cdots +Z_dS$  
has the above form $M=\mathcal E_V^\perp$ and its quotient 
$S/M$ has the form $H_V$, for a uniquely determined linear subspace 
$V\subseteq d\cdot E$.  
\end{thm}

\begin{proof} We begin with the observation that for any 
standard Hilbert module $S=G\otimes E$ 
based on a maximally symmetric graded completion $G$, 
the restriction 
of $Z_k^*$ to $\mathcal A\otimes E$ has the form 
$$
Z_k^*=u(N)\frac{\partial}{\partial z_k} , 
$$
where $u(N)$ denotes a bounded function of the number 
operator $N$ that depends on the inner product 
of $G$ (see formula (\ref{lqEq1.1})).  
It follows that for each $n\geq 0$,  the restriction 
of $Z_k^*$ to $\mathcal A_n\otimes E$ is a scalar multiple 
of $\frac{\partial}{\partial z_k}$.  

In view of (\ref{lqEq6}), 
the closure of $\mathcal E_V$ becomes  
a graded subspace of $S$ that is invariant under 
$Z_1^*, \dots,Z_d^*$.   
Thus we may conclude that 
a) $M=\mathcal E_V^\perp$ is a graded submodule of $S$, 
and b) the quotient $S/M$ is identified with 
the closure of $\mathcal E_V$, so that the 
coordinate operators 
of $S/M$ are identified with the compressions of 
$Z_1,\dots,Z_d$ to $\overline{\mathcal E_V}$.

To see that the graded 
submodule $M=\mathcal E_V^\perp$ is of degree one, note 
first that by Remark \ref{gqRem2}, 
$G\otimes V$ is a reducing submodule of $G\otimes d\cdot E=d\cdot S$, 
and Proposition \ref{lqProp1} implies that the closure of $\mathcal E_V$ 
is the space of all elements $\xi\in S=G\otimes E$ such that 
$$
(Z_1^*\xi,\dots,Z_d^*\xi)\in G\otimes V.  
$$
Let $Q\in\mathcal B(d\cdot E)$ be the projection of $d\cdot E$ onto $V^\perp$.  
If we 
realize elements of $d\cdot E$ as column vectors with components 
in $E$, then 
the preceding formula makes the assertion  
\begin{equation}\label{lqEq4}
\xi\in M^\perp \iff (\mathbf 1\otimes Q)
\begin{pmatrix}
Z_1^*\xi\\\vdots\\Z_d^*\xi
\end{pmatrix}
=0.
\end{equation}
We can realize 
$Q$ as a $d\times d$ matrix $(Q_{ij})$ of operators in $\mathcal B(E)$.  
Letting $L_i$ be the operator in $\mathcal B(S)$ defined by 
$$
L_i=Z_1^*\otimes Q_{i1}+Z_2^*\otimes Q_{i2}+\cdots+Z_d^*\otimes Q_{id}, 
\qquad i=1,\dots,d,
$$
then (\ref{lqEq4}) becomes the assertion 
$$
M^\perp = \ker L_1\cap \cdots\cap \ker L_d,   
$$
or equivalently, 
$$
M=\overline{\ran L_1^*+\cdots+\ran L_d^*}.  
$$
Finally, after noting that each operator $L_i^*$ has the 
form 
$$
L_i^*=Z_k\otimes Q_{1i}+Z_2\otimes Q_{2i}+\cdots Z_d\otimes Q_{di}, 
$$
it follows by inspection that the closure of $\ran L_i^*$ is a 
graded degree-one 
submodule of $Z_1S+\cdots+Z_dS$.  Hence these $d$ submodules generate  
a graded degree-one submodule $M\subseteq Z_1S+\cdots+Z_dS$.   

Conversely, let $M$ be a degree-one submodule of $Z_1S+\cdots+Z_dS$, 
so that $M=[\mathcal AM_1]$ where $M_1$ is a linear subspace of $S_1$.   
Consider the 
homomorphism of $\mathcal A$-modules $L: d\cdot S\to S$ defined by 
$$
L(\eta_1,\dots,\eta_d)=Z_1\eta_1+\cdots+Z_d\zeta_d.  
$$
$L$ is a bounded linear map with closed range $Z_1S+\cdots +Z_dS$.  
By normalization 
of the basis $z_1,\dots,z_d$ for $G_1$, it follows that for every 
$\xi\in M_1\subseteq G_1\otimes E$, the elements $\zeta_k\in E$ defined 
by $1\otimes \zeta_k=Z_k^*\xi$ satisfy 
$$
\xi=z_1\otimes \zeta_1+\cdots+z_d\otimes \zeta_d=L(\zeta_1,\dots,\zeta_d).  
$$
Indeed, the set $W=\{(\zeta_1,\dots,\zeta_d)\in d\cdot E:L(\zeta_1,\dots,\zeta_d)\in M_1\}$ 
is a linear 
subspace of $d\cdot E$ such that $L(W)=M_1$.  
It follows that $N=[\mathcal AW]$ is a {\em reducing}  
submodule of $d\cdot S$ with the property that the restriction 
$L_0$ of $L$ to 
$N=[\mathcal AW]$ carries $N$ onto a dense linear submanifold 
of $M$.  Note, however,  that {\em we 
cannot assert that $L(N)$ is closed}. 

In any case, we have 
$$
M^\perp=L(N)^\perp=\ker(P_NL^*)=\{\xi\in S: L^*\xi\in N^\perp\}.  
$$ 
After noting that 
$N^\perp=G\otimes W^\perp$ and that $L^*: S\to d\cdot S$ is the operator 
$L^*\xi=(Z_1^*\xi,\dots,Z_d^*\xi)$, we find that  
\begin{equation}\label{lqEq5}
M^\perp=\{\xi\in E: (Z_1^*\xi,\dots,Z_d^*\xi)\in G\otimes V\}, 
\end{equation}
where $V$ is the orthocomplement of $W$ in $d\cdot E$.  
Finally, since $M^\perp$ is the closure 
of the set of polynomials it contains, we 
have exhibited $M^\perp$ as the closure of the set of polynomials 
$f\in\mathcal A\otimes E$ with the property 
$$
D_Gf\in G\otimes V.  
$$

Now by (\ref{lqEq1.1}), $D_G$ has the form $D_G=u(N)\nabla$ where 
$u$ is a bounded function of the number operator, where $u$ 
has the further property $u(k)>0$ for every $k=0,1,2,\dots$.  Since 
$G\otimes V$ is graded and graded submodules are invariant under 
functions of the number operator, it follows that 
$$
M^\perp = \overline{\{f\in \mathcal A\otimes E: \nabla f\in \mathcal A\otimes V\}}.  
$$
Proposition \ref{lqProp1} provides a pointwise 
description of the relation $\nabla f\in\mathcal A\otimes V$
that exhibits $M^\perp$ as the closure of $\mathcal E_V$,  as 
asserted in (\ref{lqEq3}).  
\end{proof}

It is quite easy to adapt the proof of Theorem \ref{lqThm1} 
so as to give a concrete description of $G(1)$-quotients 
for arbitrary graded completions $G$.  The general statement 
is somewhat less elementary, in that one must replace the natural 
differentiation operators $\frac{\partial}{\partial z_k}$ 
with $Z_k^*$, $1\leq k\leq d$, and $\nabla$ with the operator 
$D_G$ of (\ref{lqEq0.1}).  Thus, in order to apply the more 
general result to 
a given graded completion $G$,  
one would first have to identify the operators $Z_k^*$ 
more explicitly 
as differential operators involving $\frac{\partial}{\partial z_k}$.  
We merely state the general  result here, 
leaving details for the reader.

\begin{thm}\label{lqThm2}
Let $S=G\otimes E$ be a standard Hilbert module based on 
an arbitrary graded completion $G$.  The most general 
degree $1$ submodule of $Z_1S+\cdots+Z_dS$ has the form 
$M=\mathcal E_V^\perp$, where $V$ is a linear subspace
of $d\cdot E$, $\mathcal E_V$ is defined by   
$$
\mathcal E_V=\{f\in\mathcal A\otimes E: D_Gf(z)\in V,\quad z\in\mathbb C^d\}, 
$$
and $D_G:\mathcal A\otimes E\to \mathcal A\otimes d\cdot E$ 
is the row operator $(Z_1^*,\dots,Z_d^*)$. 
Thus the corresponding 
$G(1)$-quotient is identified with the 
closure of $\mathcal E_V$ in $G\otimes E$, whose 
operators are the compressions of $Z_1,\dots,Z_d$ to $\overline{\mathcal E_V}$.  
\end{thm}

\section{Concluding Remarks, open Problems}\label{S:cp}

There is a range of conjectures associated with 
the basic problem discussed above.  The most conservative 
of them is formulated below as Conjecture A.  Since there is 
too little evidence to support the strongest conjecture  
one might entertain, we have formulated that 
as a question in Problem B.  

\vskip0.1in
{\bf Conjecture A.}  Let $G$ be a maximally symmetric 
essentially normal graded completion (such as the Hardy 
or Bergman modules of the unit ball or the space of the $d$-shift).  Then 
for every finite-dimensional Hilbert space $E$ and every 
subspace $V\subseteq d\cdot E$, the $G(1)$-quotient 
$H_V=(G\otimes E)/\mathcal E_V^\perp$ 
described in Theorem \ref{lqThm1} is essentially 
normal.  
\vskip0.1in

\vskip0.1in
{\bf Problem B.}  Fix $p\in[1,\infty)$.  Let 
$G$ be an arbitrary graded completion that 
is  $p$-essentially normal,   
let $E$ be a finite-dimensional Hilbert space, 
and let $V$ be a subspace of $d\cdot E$.  
Is the 
$G(1)$-quotient $H_V=(G\otimes E)/\mathcal E_V^\perp$ 
$p$-essentially  normal?   
\vskip0.1in

By Corollary \ref{liCor1}, an affirmative answer to Conjecture A 
would imply that all $G$-quotients are essentially 
normal in the maximally symmetric case.  A positive 
reply to Problem B would have  
more far-reaching consequences.  

\begin{rem}[Failure of Essential Normality]
At this point we should point out that many 
familiar graded completions $G$ are {\em not} essentially 
normal; consequently one cannot expect quotients of standard 
Hilbert modules based on such $G$ to be essentially 
normal.  For example, consider the Hardy or Bergman space 
$G$ of the bidisk $D\times D$, $D=\{z\in \mathbb C: |z|<1\}$.  
In either case, $G$ is a tensor product of one-dimensional 
Hilbert modules $G_0\otimes G_0$ 
where $G_0$ is the Hardy or Bergman space of the unit disk.  
Thus in both cases, the \cstar\ generated by the coordinate 
multiplications of $G$ is the tensor product $\mathcal T\otimes \mathcal T$, 
where $\mathcal T$ is the one-dimensional Toeplitz \cstar.  
Since $\mathcal T\otimes \mathcal T$ 
admits nontrivial chains of ideals such as 
$$
\mathcal K\otimes\mathcal K\subseteq 
\mathcal K\otimes \mathcal T\subseteq 
\mathcal K\otimes\mathcal T+\mathcal T\otimes\mathcal K\subseteq 
\mathcal T\otimes \mathcal T,  
$$
$\mathcal K$ denoting the compact operators  
on $G_0$, $\mathcal T\otimes \mathcal T$ is 
not commutative modulo compact operators.  
See \cite{dougHoQP}, \cite{DougId} for further discussion.  
The Bergman spaces of more general domains can give rise to 
type $I$ \cstar s with arbitrarily long composition series \cite{upAn}.  

\end{rem}

There is considerable evidence to support Conjecture A 
and its consequences.  
The results of \cite{arvpsum} imply that $H^2$-quotients 
$(H^2\otimes E)/M$ are $p$-essentially normal for $p>d$ 
whenever $M$ is generated by monomials.  Douglas 
\cite{dougPsum} has generalized that result to the context 
of more general weighted shifts.  In dimension 
$d=2$, Guo \cite{guoDef} obtained 
trace estimates which imply that $H^2/M$ is essentially 
normal for every graded submodule $M\subseteq H^2$.  
In very recent work with Wang \cite{guoWang} that result is 
improved; the new version implies that Conjecture 
A is true in dimension $d=2$.   

Finally, we have shown that Conjecture A itself 
holds in certain special 
cases (in arbitrary dimensions), 
including a) that in which $V$ is one-dimensional, 
b) that in which $V$ is of codimension $1$ in $d\cdot E$,   
and c) that in which the submodule $M=\mathcal E_V^\perp$ is 
``diagonal" in the sense that $Z_k^*M_1\perp Z_j^*M_1$ when $k\neq j$.  
Unfortunately, 
none of the three proofs appears 
to generalize to the full context of Conjecture A.  

\vskip0.1in
{\em Acknowledgement.}  I want to thank Kunyu Guo for 
pointing out an error in a lemma a previous draft 
of this paper: a condition that was asserted to hold for 
$n=1,2,\dots$ actually holds only for $n=2,3,\dots$.  
That opened a gap in the proof of the main result 
that remains unfilled at the time of this writing.  
The current version of this 
paper has been revised and reorganized 
in essential ways, and contains a more 
modest main result.

\appendix
\section{Examples with Maximum Symmetry}\label{S:ex}

We find all graded 
completions that are essentially normal and have 
maximum possible symmetry.  This means that 
the associated inner product on $\mathcal A$ 
is invariant under 
the action of the full unitary group of $\mathbb C^d$
on $\mathcal A$.  
We also show that every graded completion 
with maximum symmetry gives rise to operator $d$-tuples 
that satisfy commutation relations of a general type.  
These commutation relations are critical for the 
results in later sections.  

Examples of maximally symmetric graded 
completions include the module $H^2$ of the $d$-shift, 
the Hardy module on the unit sphere in complex dimension $d$, 
the Bergman module on 
the ball, and many others related to domains with 
symmetry that are not tied to the unit sphere.  
Of course, there is a vast array of 
more general standard Hilbert 
modules that have less symmetry, 
even examples with minimum symmetry 
in the sense that only the center of the unitary group 
acts naturally - minimum symmetry being necessary 
as part of 
the definition of graded inner product.  
However, while the problem of classifying standard 
Hilbert modules in general appears  
difficult, we are 
optimistic about further progress in 
analyzing well-chosen intermediate subclasses.

The unitary group $G=\mathcal U_d$ of $\mathbb C^d$ acts 
naturally on the algebra $\mathcal A$ of polynomials 
in $d$ variables.  One sees this most clearly by realizing 
$\mathcal A$ as the symmetric tensor algebra 
$$
\mathbb C\dotplus Z\dotplus Z^{(2)}\dotplus\cdots 
$$
over the one-particle space $Z=\mathbb C^d$, $Z^{(n)}$ denoting 
the symmetric tensor product of $n$ copies of $Z$, which 
of course can be 
identified with the space 
$\mathcal A_n\subseteq \mathbb C[z_1,\dots,z_d]$ of homogeneous polynomials 
of degree $n$ once one specifies a basis for $Z$.  The action 
of $G$ is given by second quantization
$$
\Gamma(U)=\mathbf 1_\mathbb C\dotplus U\dotplus U^2\dotplus\cdots, 
\qquad U\in G, 
$$ 
where $U^n$ denotes the restriction of 
$U^{\otimes n}\in\mathcal B(Z^{\otimes n})$ 
to the symmetric subspace $\mathcal A_n\subseteq Z^{\otimes n}$, 
$n=1,2,\dots$.  

There are many graded completions in dimension $d$ that 
are associated with rotationally-invariant measures on $\mathbb C^d$.
Indeed, let $\mu$ be a compactly supported probability measure 
on $\mathbb C^d$ that is invariant under the action of the full unitary 
group $\mathcal U_d$.  The closed subspace $G\subseteq L^2(\mu)$ 
generated by the polynomials defines a Hilbert module which, 
under appropriate mild conditions on the measure $\mu$, 
is a $p$-essentially normal graded completion for every $p>d$ 
(see Theorem \ref{exThm1} and Proposition \ref{exPropPsum}).  
The Hardy and Bergman 
modules in dimension $d$ are of this type.  

These examples 
are obviously subnormal Hilbert modules.  On the other hand, 
while the module $H^2$ of the 
$d$-shift is not subnormal and cannot be associated 
with a measure, it is also a symmetric graded completion
\cite{arvSubalgIII}.

Choose a $G$-invariant inner product $\langle\cdot,\cdot\rangle$ 
on $\mathcal A$.  Such an 
inner product  
is of course graded, so that   
$\langle \mathcal A_m,\mathcal A_n\rangle=\{0\}$ if $m\neq n$.  
Moreover, since the restriction of $\Gamma$ to 
each homogeneous subspace $\mathcal A_n$ is 
an irreducible representation of $G$, 
any two $G$-invariant inner products on 
$\mathcal A_n$ must be proportional.  Hence 
there is a sequence of positive constants 
$c_0,c_1,\dots$ such that 
$$
\langle f,g\rangle=c_n\langle f,g\rangle_{H^2},\qquad 
f,g\in \mathcal A_n,\quad n=0,1,2,\dots,   
$$
where $\langle\cdot,\cdot\rangle_{H^2}$ denotes 
the inner product of the symmetric Fock space.  
Conversely, given any sequence $c_0,c_1,\dots$ of positive 
numbers, the preceding formula defines a $G$-invariant 
inner product $\langle\cdot,\cdot\rangle$ on $\mathcal A$.  
Thus, we seek to determine all sequences $c_0,c_1,c_2,\dots$ 
with the property that the associated inner product 
leads to a standard Hilbert module that is essentially 
normal, or more generally, that is $p$-essentially normal 
for some $p>d$. 

Fixing a sequence of positive numbers $c_0,c_1,\dots$, 
it is more convenient to work with another sequence 
$\rho_0,\rho_1,\rho_2,\dots$ defined by 
$$
\rho_k=\sqrt\frac{c_{k+1}}{c_k},\qquad k=0,1,2,\dots.     
$$
Thus, $c$ and $\rho$ are related by 
$c_{k+1}=(\rho_0\rho_1\cdots\rho_k)^2c_0$,
$k=0,1,2,\dots$, so that knowing the  
sequence $\rho_0,\rho_1,\cdots$ 
is equivalent to knowing the inner 
product up to a positive scaling factor.  

\begin{thm}\label{exThm1}
Let $\rho_0,\rho_1,\dots$ be a sequence of 
positive numbers, let $G$ be the Hilbert 
space obtained by completing $\mathcal A$ 
in the symmetric inner product associated with $(\rho_k)$ 
as above.  Then the coordinate operators 
$Z_1, \dots,Z_d$ are bounded and 
$Z_1G+\cdots+Z_dG$ is closed iff 
\begin{equation}\label{exEq6}
0<\epsilon\leq \rho_k\leq M<\infty, \qquad k=0,1,2,\dots,
\end{equation}
for some positive constants $\epsilon, M$.  
%\end{enumerate}
%
In that event, 
$G$ is also essentially normal 
iff the sequence $\rho_n$ oscillates slowly 
in the sense that 
\begin{equation}\label{exEq4}
\lim_{n\to\infty} (\rho_{n+1}-\rho_n)=0.      
\end{equation}
 \end{thm}

\begin{proof}
Let $\langle \cdot,\cdot\rangle$ be the inner product 
on $\mathcal A$ associated with 
$$
c_k=(\rho_0\rho_1\cdots\rho_k)^2c_0,\qquad k=0,1,2\dots,
$$ 
and let 
$\langle\cdot,\cdot\rangle_{H^2}$ be the inner product of 
the symmetric Fock space $H^2$.  
Let $(S_1,\dots,S_d)\in\mathcal B(H^2)$ be the $d$-shift, 
let $E_n\in\mathcal B(H^2)$ be the projection onto 
the space of homogeneous polynomials $\mathcal A_n$, $n=0,1,2,\dots$, 
and let $\Delta$ be the following diagonal 
operator in $\mathcal B(H^2)$ 
\begin{equation}\label{exEq3}
\Delta=\sum_{n=0}^\infty \rho_n E_{n+1}.  
\end{equation}

We claim first that, up to a graded 
unitary equivalence, the $d$-tuple 
$(Z_1,\dots,Z_d)$ acting on $\mathcal A\subseteq G$, 
is the ``weighted $d$-shift" $(\Delta S_1,\dots,\Delta S_d)$, 
considered as a densely defined operator acting on 
$\mathcal A\subseteq H^2$.   
Indeed, we have 
\begin{equation}\label{exEq2}
\langle f,g\rangle = c_n\langle f,g\rangle_{H^2}, 
\qquad f,g\in \mathcal A_n,\quad n=0,1,2,\cdots.  
\end{equation}
Letting $W:\mathcal A\to \mathcal A$ 
be the linear map
$$
W=\sum_{n=0}^\infty \sqrt {c_n}E_n, 
$$  
one sees that 
$W$ is a linear isomorphism of $\mathcal A$ onto itself,  
and by (\ref{exEq2}) we can take $G$ to be the 
completion of $\mathcal A$ in the inner product 
$f,g\mapsto \langle Wf,Wg\rangle_{H^2}$.  
For $S_1,\dots,S_d$ and $Z_1,\dots,Z_d$ as above we 
have $Z_kf=z_k\cdot f=S_kf$ 
for polynomials $f$, hence 
$$
\langle Z_kf,g\rangle=\langle WS_kf,Wg\rangle_{H^2}=
\langle WS_kW^{-1}Wf,Wg\rangle_{H^2}
$$
and it follows that the $d$-tuple of 
restrictions of $Z_1,\dots,Z_d$ 
to $\mathcal A\subseteq H$ is
unitarily equivalent 
to the $d$-tuple of 
restrictions of $WS_1W^{-1},\dots,WS_dW^{-1}$ 
to $\mathcal A\subseteq H^2$.  
Using the commutation formula 
$S_kE_n=E_{n+1}S_k$, one 
can now compute in the obvious way to obtain 
$WS_kW^{-1}=\Delta S_k$, as asserted.  

Thus we may take $Z_k=\Delta S_k$, $k=1,\dots,d$.  
Since $S_1S_1^*+\cdots+S_dS_d^*$ is the 
projection $E_0^\perp=E_1+E_2+\cdots$, one finds 
that 
$$
Z_1Z_1^*+\cdots +Z_dZ_d^*=
\Delta (S_1S_1^*+\cdots+S_dS_d^*)\Delta=\Delta^2=\sum_{n=0}^\infty \rho_n^2E_{n+1}, 
$$
from which the equivalences characterized by (\ref{exEq6}) are apparent.  

We now consider essential normality of 
the $Z_k=\Delta S_k$, $1\leq k\leq d$.  
Noting that the commutation formula 
$S_kE_n=E_{n+1}S_k$ implies that $\Delta$ commutes with 
$S_kS_j^*$, and that  $\Delta S_k=S_k\tilde \Delta$ where 
\begin{equation}\label{exEq3a}
\tilde \Delta=\sum_{n=0}^\infty \rho_nE_n, 
\end{equation}
it follows that each commutator $[Z_j^*,Z_k]$ can be 
written 
\begin{equation*}
S_j^*\Delta^2S_k-\Delta S_kS_j^*\Delta=S_j^*S_k\tilde \Delta^2-S_kS_j^*\Delta^2=
[S_j^*,S_k]\tilde \Delta^2+S_kS_j^*(\tilde \Delta^2-\Delta^2).  
\end{equation*}
Since the self-commutators 
$[S_j^*,S_k]$ are known to belong to $\mathcal L^p$ 
for every $p>d$ \cite{arvSubalgIII} and $\tilde \Delta^2$ is bounded, 
$[S_j^*,S_k]\tilde \Delta^2$ is compact.  Hence 
$[Z_j^*,Z_k]$ is compact for all $j,k$ iff 
$S_kS_j^*(\tilde \Delta^2-\Delta^2)$ is compact for 
all $k,j$.  Noting again that $S_1S_1^*+\cdots+S_dS_d^*$ is 
a rank-one perturbation of the identity, we find  that 
$S_kS_j^*(\tilde \Delta^2-\Delta^2)$ is compact for all $j,k$ iff 
$\tilde \Delta^2-\Delta^2$ is compact.  Since 
$$
\tilde \Delta^2-\Delta^2=(\tilde \Delta+\Delta)
(\tilde \Delta-\Delta)=(\tilde \Delta+\Delta)
(\rho_0E_0+\sum_{n=1}^\infty(\rho_n-\rho_{n-1})E_n)
$$
and $E_0,E_1,E_2,\dots$ is a sequence of 
mutually orthogonal finite-dimensional 
projections, (\ref{exEq4}) follows. 
\end{proof}

The 
examples characterized in Theorem \ref{exThm1} 
include $H^2$ (take $\rho_k=1$ for every $k$), 
the Hardy module, which is associated with the 
sequence
$$
\rho_k = \sqrt\frac{k+1}{k+d},\qquad k=0,1,2,\dots,   
$$
and the Bergman module.  
While these examples are 
all associated with the unit sphere in the sense 
that $Z_1Z_1^*+\cdots+Z_dZ_d^*$ is a compact 
perturbation of the identity $\mathbf 1$,  
there are many others that are not.  

For example, for 
any two   
positive constants $0<r_1<r_2<\infty$, the slowly oscillating 
sequence $\rho_0,\rho_1,\dots$ defined by 
\begin{equation}\label{exEq1}
\rho_k^2 = r_1 + (r_2-r_1)\frac{1+\sin \sqrt k}{2},\qquad k=0,1,2, \dots 
\end{equation}
defines an essentially normal standard Hilbert module 
with the property that the  
spectrum of $Z_1Z_1^*+\cdots+Z_dZ_d^*$ is the union 
$\{0\}\cup [r_1,r_2]$.  Such a Hilbert module is associated 
with the annular region in $\mathbb C^d$ 
$$
X=\{z\in \mathbb C^d: r_1\leq \|z\|\leq r_2\},
$$  
since it gives rise to an exact sequence of \cstar s
$$
0\longrightarrow \mathcal K\longrightarrow C^*(Z_1,\dots,Z_d)
\longrightarrow C(X)\longrightarrow 0.  
$$

\begin{rem}[$p$-essential normality]
Using the fact that the Hilbert module 
$H^2$ is $p$-essentially normal for every 
$p>d$ and that the dimension of $\mathcal A_n$ grows 
as a polynomial of degree $d-1$, 
is not hard to adapt the proof of Theorem \ref{exThm1} 
to establish the following characterization 
of $p$-essential normality:

\begin{prop}\label{exPropPsum}
Let $\rho_0,\rho_1,\dots$ be a 
sequence satisfying (\ref{exEq4}) and 
property (ii) of Theorem \ref{exThm1}.  Let $G$ 
be the graded completion obtained from the associated 
inner product.  For every $p\in(d,\infty)$, the following 
are equivalent:
\begin{enumerate}
\item[(i)]$G$ is $p$-essentially normal.  
\item[(ii)]
\begin{equation}\label{exEq5}
\sum_{k=1}^\infty k^{d-1}|\rho_{k+1}-\rho_k|^p <\infty.  
\end{equation}
\end{enumerate}
\end{prop}

We omit the proof since we do not require this result.  

However, we point out that given some number $\alpha>d$, it 
is easy to use Proposition \ref{exPropPsum} to 
find examples of sequences $(\rho_k)$ 
that give rise to graded completions $G$ that 
are $p$-essentially normal for all $p>\alpha>d$ but 
not for $p$ in the range $d<p\leq \alpha$.  For instance, 
since the sequence of (\ref{exEq1})
behaves so that $|\rho_{k+1}-\rho_k|=O(k^{-1/2})$, 
straightforward estimates using (\ref{exEq5}) 
show that the Hilbert module $S$ 
associated with that sequence is  
$p$-essentially normal iff $p> 2d$.  In particular,  
this $S$ is {\em not} $p$-essentially normal for $p$ 
in the range $d<p\leq 2d$.  
\end{rem}

%\vfill

%\bibliography{bibData}  %Remove this line when finished, see below.  

\begin{thebibliography}{DMV00}

\bibitem[Arv98]{arvSubalgIII}
W.~Arveson.
\newblock Subalgebras of {$C^*$}-algebras {III}: Multivariable operator theory.
\newblock {\em Acta Math.}, 181:159--228, 1998.
\newblock arXiv:funct-an/9705007.

\bibitem[Arv00]{arvCurv}
W.~Arveson.
\newblock The curvature invariant of a {H}ilbert module over {${\mathbb
  C}[z_1,\dots,z_d]$}.
\newblock {\em J. Reine Angew. Mat.}, 522:173--236, 2000.
\newblock ar{X}iv:math.OA/9808100.

\bibitem[Arv02]{arvDirac}
W.~Arveson.
\newblock The {D}irac operator of a commuting $d$-tuple.
\newblock {\em Jour. Funct. Anal.}, 189:53--79, 2002.
\newblock ar{X}iv:math.OA/0005285.

\bibitem[Arv04]{arvpsum}
W.~Arveson.
\newblock $p$-summable commutators in dimension $d$.
\newblock {\em J. Oper. Th.}, 2004.
\newblock ar{X}iv:math.OA/0308104 v2.

\bibitem[BDF77]{bdf}
L.G. Brown, R.~Douglas, and P.~Fillmore.
\newblock Extensions of {$C^*$}-algebras and {$K$}-{H}omology.
\newblock {\em Ann Math.}, 105(2):265--324, March 1977.

\bibitem[Cur81]{curFred}
R.~Curto.
\newblock Fredholm and invertible {$n$}-tuples of operators. the deformation
  problem.
\newblock {\em Trans. A.M.S.}, 266(1):129--159, 1981.

\bibitem[DH71]{dougHoQP}
R.~G. Douglas and Roger Howe.
\newblock On the {$C^*$}-algebra of {T}oeplitz operators on the quarterplane.
\newblock {\em Trans. A.M.S.}, 158:203--217, 1971.

\bibitem[DM03]{dougMisEQ}
R.~G. Douglas and G.~Misra.
\newblock Equivalence of quotient {H}ilbert modules.
\newblock {\em Proc. Indian Acad. Sci. (Math. Sci.)}, 113(3):281--291, August
  2003.

\bibitem[DMV00]{dougMV1}
R.~G. Douglas, G.~Misra, and C.~Varughese.
\newblock On quotient modules - the case of arbitrary multiplicity.
\newblock {\em J. Funct. Anal.}, 174:364--398, 2000.

\bibitem[Dou05a]{dougPsum}
R.~G. Douglas.
\newblock Essentially reductive {H}ilbert modules.
\newblock {\em J. Oper. Th.}, 2005.
\newblock arXiv:math.OA/0404167.

\bibitem[Dou05b]{DougId}
R.~G. Douglas.
\newblock Ideals in {T}oeplitz algebras.
\newblock {\em Houston Math. J.}, 2005.
\newblock to appear.

\bibitem[GS05]{GRS}
S.~Gleason, J.~Richter and C.~Sundberg.
\newblock On the index of invariant subspaces in spaces of analytic functions
  in several complex variables.
\newblock {\em to appear in Crelle's Journal}, 2005.

\bibitem[Guo03]{guoDef}
K.~Guo.
\newblock Defect operators for submodules of {$H^2_d$}.
\newblock {\em preprint}, 2003.

\bibitem[GW05]{guoWang}
K.~Guo and K.~Wang.
\newblock Essentially normal {H}ilbert modules and {$K$}-homology.
\newblock {\em preprint}, 2005.

\bibitem[MV93]{mulVas}
V.~M{\"u}ller and F.-H. Vasilescu.
\newblock Standard models for some commuting multioperators.
\newblock {\em Proc. Amer. Math. Soc.}, 117:979--989, 1993.

\bibitem[Tay70a]{taylor1}
J.~L. Taylor.
\newblock The analytic functional calculus for several commuting operators.
\newblock {\em Acta Math}, 125:1--38, 1970.
\newblock MR 42 6622.

\bibitem[Tay70b]{taylor2}
J.~L. Taylor.
\newblock A joint spectrum for several commuting operators.
\newblock {\em J. Funct. Anal.}, 6:172--191, 1970.
\newblock MR 42 3603.

\bibitem[Upm84]{upAn}
Harald Upmeier.
\newblock Toeplitz {$C^*$}-algebras on bounded symmetric domains.
\newblock {\em Ann. Math.}, 119(3):549--576, 1984.

\end{thebibliography}
\bibliographystyle{alpha}
\newcommand{\noopsort}[1]{} \newcommand{\printfirst}[2]{#1}
  \newcommand{\singleletter}[1]{#1} \newcommand{\switchargs}[2]{#2#1}

\end{document}